\documentclass{article}
\usepackage{graphicx}
\usepackage{amsmath}
\usepackage{amsfonts}
\usepackage{amssymb}
\newtheorem{theorem}{Theorem}

\newtheorem{corollary}[theorem]{Corollary}

\newtheorem{definition}[theorem]{Definition}

\newtheorem{lemma}[theorem]{Lemma}

\newtheorem{proposition}[theorem]{Proposition}

\begin{document}

\title{Lagrangian submanifolds in \\Hyperk\"{a}hler manifolds,\\Legendre transformation}
\author{Naichung Conan Leung}
\maketitle
\begin{abstract}
We develop the foundation of the \textit{complex symplectic geometry} of
Lagrangian subvarieties in a hyperk\"{a}hler manifold. We establish a
characterization, a Chern number inequality, topological and geometrical
properties of Lagrangian submanifolds. We discuss a category of Lagrangian
subvarieties and its relationship with the theory of Lagrangian intersection.

We also introduce and study extensively a normalized \textit{Legendre
transformation} of Lagrangian subvarieties under a birational transformation
of projective hyperkahler manifolds. We give a \textit{Pl\"{u}cker type
formula} for Lagrangian intersections under this transformation.
\end{abstract}

\pagebreak 

\bigskip

\section{Introduction}

A Riemannian manifold $M$ of real dimension $4n$ is \textit{hyperk\"{a}hler}
if its holonomy group is $Sp\left(  n\right)  $. Its has three complex
structures $I,J$ and $K$ satisfying the Hamilton relation $I^{2}=J^{2}%
=K^{2}=IJK=-1$. We will fix one complex structure $J$ on $M$%
.\footnote{Throughout this paper $M$ is a complex manifold with a fixed
complex structure $J$ unless specified otherwise and all submanifolds are
complex submanifolds with respect to $J$.} We denote its K\"{a}hler form as
$\omega$ and its holomorphic two form as $\Omega$ which is nondegenerate and
defines a (holomorphic) symplectic structure on $M$.

A submanifold $C$ in $M$ of dimension $n$ is called a \textit{Lagrangian }if
the restriction of $\Omega$ to it is zero. In the \textit{real} symplectic
geometry, Lagrangian submanifolds plays a very key role, for example in
geometric quantization, Floer theory, Kontsevich's homological mirror
conjecture and Strominger, Yau and Zaslow geometric mirror conjecture.

The objective of this article is twofold: (1) We develop the foundation of the
\textit{complex symplectic geometry} of Lagrangian submanifolds in a
hyperk\"{a}hler manifold.\footnote{Many results presented here can be applied
to any holomorphic symplectic manifold.} This subject was previously studied
by Donagi and Markman \cite{DM}, \cite{DM2}, Hitchin \cite{Hi2} and others.
(2) We introduce a Legendre transformation of Lagrangian subvarieties along
$\mathbb{P}^{n}$ in $M$ and we establish a \textit{Pl\"{u}cker type formula}.

Here we summarize our results in this paper: First it is not difficult to
recognize a Lagrangian: (i) $C$ is Lagrangian if and only if its Poincar\'{e}
dual $\left[  C\right]  $ represents a $\Omega$-primitive class in
$H^{2n}\left(  M,\mathbb{Z}\right)  $. In particular being Lagrangian is
invariant under deformation. This also enable us to define singular Lagrangian
subvarieties. (ii) When $C$ is a complete intersection in $M$, we have
\[
\int_{C}c_{2}\left(  M\right)  \omega^{n-2}\geq0\text{,}%
\]
moreover the equally sign holds if and only if $C$ is a Lagrangian. (iii) the
rational cobordism class of the manifold $C$ is completely determined by the
cohomology class $\left[  C\right]  \in H^{2n}\left(  M,\mathbb{Z}\right)  $.

\bigskip

Second the second fundamental form of $C$ in $M$ defines a cubic vector field
$\widetilde{A}\in\Gamma\left(  C,Sym^{3}T_{C}\right)  $ and a cubic form,
\begin{align*}
\mathbf{c}_{C} &  :Sym^{3}H^{0}\left(  C,T_{C}^{\ast}\right)  \rightarrow
\mathbb{C}\\
\mathbf{c}_{C}\left(  \phi,\eta,\zeta\right)   &  =\int_{C}\phi_{i}\eta
_{j}\zeta_{k}\widetilde{A}^{ijk}\frac{\omega^{n}}{n!}\text{.}%
\end{align*}
$H^{0}\left(  C,T_{C}^{\ast}\right)  $ can be identified as the tangent space
of the moduli space $\mathcal{M}$ of Lagrangian submanifolds. By varying $C$,
the cubic tensor $\mathbf{c}_{C}$ gives a holomorphic cubic tensor
$\mathbf{c}\in H^{0}\left(  \mathcal{M},Sym^{3}T_{\mathcal{M}}^{\ast}\right)
$. This cubic form defines a torsion free flat symplectic connection $\nabla$
on the tangent bundle of $\mathcal{M}$ and it satisfies $\nabla\wedge
J_{\mathcal{M}}=0$. Namely it is a special K\"{a}hler structure on
$\mathcal{M}$.

\bigskip

Third the category $\mathcal{C}_{M}$ of Lagrangian subvarieties in $M$ with
$Hom_{\mathcal{C}_{M}}\left(  C_{1},C_{2}\right)  =\Sigma\left(  -1\right)
^{q}Ext_{O_{M}}^{q}\left(  O_{C_{1}},O_{C_{2}}\right)  $ refines the
intersection theory of Lagrangians. For example
\[
\dim Hom_{\mathcal{C}_{M}}\left(  C_{1},C_{2}\right)  =\left(  -1\right)
^{n}C_{1}\cdot C_{2}\text{.}%
\]
When $C_{1}$ and $C_{2}$ intersects cleanly, this equals the Euler
characteristic of the intersection up to sign,
\[
C_{1}\cdot C_{2}=\pm e\left(  C_{1}\cap C_{2}\right)  \text{.}%
\]
For each individual Ext group, we have (i) when $C$ is a complete intersection
Lagrangian submanifold in $M$ then
\[
Ext_{O_{M}}^{q}\left(  O_{C},O_{C}\right)  =H^{q}\left(  C,\mathbb{C}\right)
\text{, for all }q\text{,}%
\]
(ii) when $C_{1}$ and $C_{2}$ intersects transversely then $Ext_{O_{M}}%
^{q}\left(  O_{C_{1}},O_{C_{2}}\right)  =0$ for $q<n$ and equals
$\mathbb{C}^{s}$ for $q=n$.

We also study the derived category of Lagrangian coherent sheaves on $M$,
denote $D_{Lag}^{b}\left(  M\right)  $. For example the structure sheaf of any
Lagrangian subvariety in $M$ defines an object in this category.

\bigskip

Fourth we have good understanding of coisotropic subvarieties and their
corresponding symplectic reductions in hyperk\"{a}hler manifolds. Coisotropic
submanifolds share some of Lagrangian properties. They can be characterized by
the $\Omega$-primitivity property of their Poincar\'{e} duals. There is also a
Chern number inequality in the complete intersection case. Furthermore for a
generic complex structure in the twistor family of $M $ there is no isotropic
or coisotropic subvarieties.

In \textit{real} symplectic geometry, any well-behaved coisotropic submanifold
gives rise to a reduction $\pi_{D}:D\rightarrow B$ with $B$ another symplectic
manifold. Lagrangians in $M$ can be reduced to Lagrangians in $B$ or projected
to Lagrangians of $M$ inside $D$. In the \textit{complex }case we can define a
reduction functor $\mathbf{R}_{D}:D_{Lag}^{b}\left(  M\right)  \rightarrow
D_{Lag}^{b}\left(  B\right)  $ and a projection functor $\mathbf{P}%
_{D}:D_{Lag}^{b}\left(  M\right)  \rightarrow D_{Lag}^{b}\left(  M\right)  $
using Fourier-Mukai type functors.

Suppose the reduction process occurs inside $M$, namely there is a birational
contraction $\pi:M\rightarrow Z$ such that $D$ is the exceptional locus and
$B$ is the discriminant locus. In this case $\pi_{D}:D\rightarrow B $ is a
$\mathbb{P}^{k}$-bundle and its relative cotangent bundle is the normal bundle
of $D$ in $M$ (\cite{Wi}, \cite{Na}, \cite{HY}). Moreover one can replace $D$
by its dual $\mathbb{P}^{k}$-bundle and produce another holomorphic symplectic
manifold $M^{\prime}$, called the Mukai elementary modification \cite{Mu1}%
.\footnote{There is also stratified version of this construction by Markman
\cite{Mar}.}

\bigskip

Fifth we can define Legendre transformation on any Mukai elementary
modification. For example when $D\cong\mathbb{P}^{n}$ we associate to each
Lagrangian subvariety $C$ in $M$ (not equals to $D$) another Lagrangian
subvariety $C^{\vee}$ in $M^{\prime}$.

We will explain its relationship with the classical Legendre transformation.
Roughly speaking it is the hyperk\"{a}hler quotient by $S^{1}$ of the Legendre
transformation of defining functions of $C$ on the linear symplectic manifold
$T^{\ast}\mathbb{C}^{n+1}$. This can also be regarded as a generalization of
the dual varieties construction.

We establish the following \textbf{Pl\"{u}cker type formula }for the Legendre
transformation,
\[
C_{1}\cdot C_{2}+\frac{\left(  C_{1}\cdot\mathbb{P}^{n}\right)  \left(
C_{2}\cdot\mathbb{P}^{n}\right)  }{\left(  -1\right)  ^{n+1}\left(
n+1\right)  }=C_{1}^{\vee}\cdot C_{2}^{\vee}+\frac{\left(  C_{1}^{\vee}%
\cdot\mathbb{P}^{n\ast}\right)  \left(  C_{2}^{\vee}\cdot\mathbb{P}^{n\ast
}\right)  }{\left(  -1\right)  ^{n+1}\left(  n+1\right)  }\text{.}%
\]
When $n=2$ this formula is essentially equivalent to a classical Pl\"{u}cker
formula for plane curves.

Motivated from this formula we define a \textit{normalized Legendre
transformation}
\begin{align*}
\mathcal{L}\left(  C\right)   &  =C^{\vee}+\frac{\left(  C\cdot\mathbb{P}%
^{n}\right)  +\left(  -1\right)  ^{n+1}\left(  C^{\vee}\cdot\mathbb{P}^{n\ast
}\right)  }{n+1}\mathbb{P}^{n\ast}\text{ if }C\neq\mathbb{P}^{n}\\
\mathcal{L}\left(  \mathbb{P}^{n}\right)   &  =\left(  -1\right)
^{n}\mathbb{P}^{n\ast}\text{.}%
\end{align*}
This transformation $\mathcal{L}$ preserves the intersection product:
\[
C_{1}\cdot C_{2}=\mathcal{L}\left(  C_{2}\right)  \cdot\mathcal{L}\left(
C_{2}\right)  .
\]
When $n=1$ we simply have $M=M^{\prime}$ and $C=C^{\vee},$ however the
normalized Legendre transformation is interesting, it coincides with the Dehn
twist along a $\left(  -2\right)  $-curve in $M$.

In the general case, the Plucker type formula and the definition of the
normalized Legendre transformation will involve the reduction and the
projection of a Lagrangian with respect to the coisotropic exceptional
submanifold $D$.

Next we are going to look at an \textit{explicit} example of a hyperk\"{a}hler
manifold and its flop.

\bigskip

\textbf{The cotangent bundle of }$\mathbb{P}^{n}$

In \cite{Ca}, \cite{Ca2} Calabi showed that the cotangent bundle of the
complex projective space is a hyperk\"{a}hler manifold and its metric can be
described explicitly as follow: Let $z^{1},...,z^{n}$ be a local inhomogeneous
coordinate system in $\mathbb{P}^{n}$ and $\zeta_{1},...,\zeta_{n}$'s be the
corresponding coordinate system on the fibers of $T^{\ast}\mathbb{P}^{n}$,
i.e. $\Sigma\zeta_{j}dz^{j}$ represents a point in $T^{\ast}\mathbb{P}^{n}$.
The symplectic form on $T^{\ast}\mathbb{P}^{n}$ is given by $\Sigma
dz^{j}\wedge d\zeta_{j}$. The hyperk\"{a}hler K\"{a}hler form on $T^{\ast
}\mathbb{P}^{n}$ is given by
\[
\omega=\partial\bar{\partial}\left(  \log\left(  1+\left|  z\right|
^{2}\right)  +f\left(  t\right)  \right)  ,
\]
where $f\left(  t\right)  =\sqrt{1+4t}-\log\left(  1+\sqrt{1+4t}\right)  $ and
$t=\left(  1+\left|  z\right|  ^{2}\right)  \left(  \left|  \zeta\right|
^{2}+\left|  z\cdot\zeta\right|  ^{2}\right)  .$ Calabi's approach is to look
for a $U\left(  n+1\right)  $-invariant hyperk\"{a}hler metric and reduces the
problem to solving an ODE for $f\left(  t\right)  $.

\bigskip

Another approach by Hitchin \cite{Hi} is to construct the hyperk\"{a}hler
structure on $T^{\ast}\mathbb{P}^{n}$ by the method of hyperk\"{a}hler
quotient, which is analogous to the symplectic quotient (or symplectic
reduction) construction. Consider $V=\mathbb{C}^{n+1}$ with the diagonal
$S^{1}$-action by multiplication. Its induced action on $\mathbb{H}%
^{n+1}=V\times V^{\ast}=T^{\ast}V$ is given by
\begin{align*}
S^{1}\times T^{\ast}V &  \rightarrow T^{\ast}V,\\
e^{i\theta}\cdot\left(  x,\xi\right)   &  =\left(  e^{i\theta}x,e^{-i\theta
}\xi\right)  \text{.}%
\end{align*}
This $S^{1}$-action preserves both the K\"{a}hler form $dx\wedge d\bar{x}%
+d\xi\wedge d\bar{\xi}$ and the natural holomorphic symplectic form $dx\wedge
d\xi$ on $T^{\ast}V$. Namely it preserves the hyperk\"{a}hler structure on
$\mathbb{H}^{n+1}=T^{\ast}V$. The real and complex moment maps are given
respectively by
\begin{align*}
\mu_{J} &  =i\left|  x\right|  ^{2}-i\left|  \xi\right|  ^{2}\in
i\mathbb{R},\\
\mu_{c} &  =\xi\left(  x\right)  \in\mathbb{C}\text{.}%
\end{align*}
We can also combine them to form the hyperk\"{a}hler moment map:
\begin{align*}
\mu &  :\mathbb{H}^{n+1}\rightarrow i\mathbb{R}+\mathbb{C}=i\mathbb{R}%
^{3}\text{.}\\
\mu &  =\left(  \mu_{J},\mu_{c}\right)  =\left(  \mu_{J},\mu_{I},\mu
_{K}\right)  .
\end{align*}

If we take $\lambda=i\left(  1,0,0\right)  \in i\mathbb{R}^{3}$ for the
hyperk\"{a}hler quotient construction, then $S^{1}$ acts freely on $\mu
^{-1}\left(  \lambda\right)  $ and the quotient is a smooth hyperk\"{a}hler
manifold $M$,
\begin{align*}
M  &  =T^{\ast}V/_{HK}S^{1}=\mu^{-1}\left(  \lambda\right)  /S^{1}\\
&  =\left\{  \left(  x,\xi\right)  :\xi\left(  x\right)  =0,\left|  x\right|
^{2}-\left|  \xi\right|  ^{2}=1\right\}  /S^{1}.
\end{align*}

We can identify $M$ with $T^{\ast}\mathbb{P}^{n}$ explicitly as follows: It is
not difficult to see that
\begin{align*}
p &  :M\rightarrow\mathbb{P}^{n}\\
\left(  x,\xi\right)   &  \rightarrow y=x\left(  1+\left|  \xi\right|
^{2}\right)  ^{-1/2}%
\end{align*}
defines a map from $M$ to $\mathbb{P}^{n}$ with a section given by $\xi=0.$
The fiber of $p$ over the point $y=\left[  1,0,\cdots,0\right]  \in
\mathbb{P}^{n}$ consists of those $\left(  x,\xi\right)  $'s of the form
\begin{align*}
x &  =\left(  a,0,\cdots,0\right) \\
\xi &  =\left(  0,b_{1},...,b_{n}\right)  \text{,}%
\end{align*}
and satisfying $\left|  a\right|  ^{2}=1+\Sigma\left|  b_{j}\right|  ^{2}$,
i.e. $p^{-1}\left(  y\right)  $ is parametrized by $\left(  b_{1}%
,...,b_{n}\right)  \in\mathbb{C}^{n}$. Note that $p^{-1}\left(  y\right)  $
can be naturally identified with $T_{y}^{\ast}\mathbb{P}^{n}=Hom\left(
T_{y}\mathbb{P}^{n},\mathbb{C}\right)  $ via
\[
\left(  0,c_{1},\cdots,c_{n}\right)  \rightarrow\Sigma b_{j}c_{j}.
\]
By using the $U\left(  n+1\right)  $-symmetry, this gives a natural
identification between $M$ and $T^{\ast}\mathbb{P}^{n}$. Moreover $\left\{
\xi=0\right\}  $ in $M$ corresponds to the zero section in $T^{\ast}%
\mathbb{P}^{n}$, we simply denote it as $\mathbb{P}^{n}$.

The only compact Lagrangian submanifold in $T^{\ast}\mathbb{P}^{n}$ is
$\mathbb{P}^{n}$. However there are many non-compact Lagrangian submanifolds.
For any submanifold $S$ in $\mathbb{P}^{n}$, its conormal bundle
$N_{S/\mathbb{P}^{n}}^{\ast}$ is a Lagrangian submanifold in $T^{\ast
}\mathbb{P}^{n}$, a well-known construction in symplectic geometry.

\bigskip

There is a natural isomorphism $T^{\ast}V\cong T^{\ast}\left(  V^{\ast
}\right)  $ because of $\left(  V^{\ast}\right)  ^{\ast}\cong V$. We can
therefore carry out the above construction on $T^{\ast}\left(  V^{\ast
}\right)  $ exactly as before and obtain another hyperk\"{a}hler manifold
$M^{\prime}$, which is of course isomorphic to $T^{\ast}\mathbb{P}^{n\ast}$.
In terms of the above coordinate system on $V$, the only difference is to
replace the original $S^{1}$-action by its inverse action. This would change
the sign of the moment maps and therefore we can identify $M^{\prime}$ as
follow,
\begin{align*}
M^{\prime} &  =T^{\ast}\left(  V^{\ast}\right)  /_{HK}S^{1}\\
&  =\left\{  \left(  \xi,x\right)  \in V^{\ast}\times V:\xi\left(  x\right)
=1,\left|  \xi\right|  ^{2}-\left|  x\right|  ^{2}=1\right\}  /S^{1}.
\end{align*}
The two holomorphic symplectic manifolds $M$ and $M^{\prime}$ are isomorphic
outside their zero sections, the birational map is given explicitly as
follow,
\begin{align*}
\Phi &  :M\dashrightarrow M^{\prime}\\
\left(  x,\xi\right)   &  \rightarrow\left(  \left|  \frac{x}{\xi}\right|
\xi,\left|  \frac{\xi}{x}\right|  x\right)  .
\end{align*}

In fact the zero section $\mathbb{P}^{n}$ inside $M=T^{\ast}\mathbb{P}^{n}$
can be blown down and we obtain a variety $M_{0}$. Both $M$ and $M^{\prime}$
are two different crepant resolutions \cite{Ca2} of the isolated singularity
of $M_{0}$ and is usually called a \textit{flop. }This is also the basic
structure in the Mukai's elementary modification \cite{Mu1}. Explicitly we
have,
\begin{align*}
\pi &  :M\rightarrow M_{0},\\
\pi\left(  \left[  x,\xi\right]  \right)   &  \rightarrow x\otimes\xi,
\end{align*}
where $M_{0}=\left\{  A\in End\left(  \mathbb{C}^{n+1}\right)  :TrA=0,\text{
rank}\left(  A\right)  \leq1\right\}  $. It is rather easy to check that $\pi$
is the blown down morphism of $\mathbb{P}^{n}$ inside $M$. The situation for
$M^{\prime}$ is identical.

\section{Isotropic and coisotropic submanifolds}

The most natural class of submanifolds in a symplectic manifold $M$ consists
of those $C$ in $M$ with the property that the restriction of the symplectic
form $\Omega$ to $C$ is as degenerate as possible. Such a submanifold $C$ is
called isotropic, coisotropic or Lagrangian according to $\dim C\leq n$, $\dim
C\geq n$ or $\dim C=n$ respectively.

\subsection{\label{Sec Iso-Coiso property}Definitions and properties}

We first look at the linear case, i.e. $M$ and $C$ are vector space and its
linear subspace respectively. The complement of $C$ in $M$ is defined as follow,%

\[
C^{\bot}=\left\{  v\in M:\Omega\left(  v,w\right)  =0\text{ for all }w\in
C\right\}  .
\]
$C$ is called \textit{isotropic }(resp. \textit{coisotropic }or
\textit{Lagrangian}) if $C\subset C^{\bot}$ (resp. $C^{\bot}\subset C$ or
$C=C^{\bot} $). Here is the standard example: when $M=\mathbb{C}^{2n}$ with
$\Omega=dz^{1}\wedge dz^{n+1}+\cdots+dz^{n}\wedge dz^{2n}$, then the linear
span of $\frac{\partial}{\partial z^{1}},\cdots,\frac{\partial}{\partial
z^{m}}$ is an isotropic (resp. coisotropic or Lagrangian) subspace if $m\leq
n$ (resp. $m\geq n$ or $m=n$). A useful linear algebra fact is any isotropic
or coisotropic subspace of $M$ is equivalent to the one in the standard
example up to an automorphism of $M$ which preserves $\Omega$, namely a
symplectomorphism. For general symplectic manifolds we have the following
standard definition.

\begin{definition}
If $\left(  M,\Omega\right)  $ is a symplectic manifold and $C$ is a
submanifold of $M$, then $C$ is called isotropic (resp. coisotropic or
Lagrangian) if $TC\subset TC^{\bot}$ (resp. $TC\supset TC^{\bot}$ or
$TC=TC^{\bot}$). Here
\[
TC^{\bot}=\left\{  v\in TM|_{C}:\Omega\left(  v,w\right)  =0\text{ for all
}w\in TC\right\}  .
\]
\end{definition}

When $\dim C=1$ (resp. $2n-1$), $C$ is automatically isotropic (resp.
coisotropic). From the definition it is clear that $C$ being isotropic is
equivalent to $\Omega|_{C}=0$. This happens if $C$ has no nontrivial
holomorphic two form, i.e. $H^{2,0}\left(  C\right)  =0$. In particular any
submanifold $C$ in $M$ with dimension greater than $n$ has $H^{2,0}\left(
C\right)  \neq0$. More generally if $\dim C=n+k$ then $H^{2j,0}\left(
C\right)  \neq0$ for $j=0,\cdots,k$.

\begin{lemma}
If $M$ is a hyperk\"{a}hler manifold and $C$ is a submanifold in $M$ of
dimension $n+k$ then the restriction of $\left(  \Omega\right)  ^{k}$ to $C$
is always nowhere vanishing and moreover $\left(  \Omega\right)  ^{k+1}=0$ if
and only if $C$ is coisotropic.
\end{lemma}

Proof of lemma: In the standard example the above assertion can be verified
directly. For the general case the assertion follows from the fact that every
coisotropic subspace in a symplectic vector space can be conjugated by a
symplectomorphism to the one in the standard example. $\blacksquare$

With the help of the above lemma we can prove the following useful result.

\begin{theorem}
If $M$ is a compact hyperk\"{a}hler manifold and $C$ is a submanifold in $M$
then
\[%
\begin{array}
[c]{ccc}%
\text{(i) } & C\text{ isotropic} & \text{iff }L_{\Omega}\left[  C\right]
=0,\\
\text{(ii) } & \text{ }C\text{ coisotropic} & \text{iff }\Lambda_{\Omega
}\left[  C\right]  =0.
\end{array}
\]
Here $\left[  C\right]  \in H^{\ast}\left(  M,\mathbb{Z}\right)  $ denotes the
Poincar\'{e} dual of $C$.
\end{theorem}

Proof: We need to use the Hard Lefschetz $sl_{2}$-action induced by
$L_{\Omega}$ and $\Lambda_{\Omega}$ (see appendix for details). For (i) when
$C$ is isotropic we have $\Omega|_{C}=0$ and therefore $L_{\Omega}\left[
C\right]  =0$ by Poincar\'{e} duality. For the converse we suppose that
$L_{\Omega}\left[  C\right]  =0$, i.e. $\left[  C\right]  \cup\Omega=0$. By
the Hard Lefschetz $sl_{2}$-action the dimension of $C$ must be strictly less
than $n$, say $m$.

If $m=1$ then $C$ is already isotropic. So we assume $m\geq2$ and we have
\[
0=\int_{M}\left[  C\right]  \Omega\bar{\Omega}\omega^{m-2},
\]
where $\omega$ is the K\"{a}hler form on $M$. This is the same as $\int
_{C}\Omega\bar{\Omega}\omega^{n-2}=0$. Recall that $\omega|_{C}$ defines a
K\"{a}hler structure on $C$. By the Riemann bilinear relation on the
K\"{a}hler manifold $C$, the function $\left(  \Omega\bar{\Omega}\omega
^{n-2}\right)  /\omega^{n}$ is proportional to the norm square of $\Omega
|_{C}$. Therefore the holomorphic form $\Omega|_{C}$ must vanish, namely $C$
is isotropic.

For (ii) we first suppose that $\Lambda_{\Omega}\left[  C\right]  =0$. By the
Hard Lefschetz $sl_{2}$-action induced by $L_{\Omega}$ and $\Lambda_{\Omega}$,
we have $\dim C=m=n+k>n$ and
\[
\left[  C\right]  \cup\left(  \Omega\right)  ^{k+1}=0\in H^{2n+2k+2}\left(
M,\mathbb{C}\right)  \text{.}%
\]
Therefore
\begin{align*}
0 &  =\int_{M}\left[  C\right]  \Omega^{k+1}\bar{\Omega}^{k+1}\omega^{2n-2}\\
&  =\int_{C}\Omega^{k+1}\bar{\Omega}^{k+1}\omega^{2n-2}.
\end{align*}
Using the Riemann bilinear relation as before, we have $\left(  \Omega\right)
^{k+1}|_{C}=0$.\ By the above lemma, $C$ is a coisotropic submanifold. The
converse is clear. Hence the result. $\blacksquare$

\bigskip

The theorem says that a submanifold of $M$ being isotropic (or coisotropic)
can be detected cohomologically. An immediate corollary is that such property
is invariant under deformation. It also enables us to define isotropic (or
coisotropic) subvarieties which might be singular, or even with non-reduced
scheme structure.

\begin{definition}
Suppose $M$ is a compact hyperk\"{a}hler manifold. A subscheme $C$ of $M$ is
called isotropic (resp. coisotropic) if $L_{\Omega}\left[  C\right]  =0$
(resp. $\Lambda_{\Omega}\left[  C\right]  =0$).
\end{definition}

The next proposition gives a good justification of this definition.

\begin{proposition}
Suppose $M$ is a compact hyperk\"{a}hler manifold and $C$ is an irreducible
subvariety. Then $C$ is isotropic (resp. coisotropic) if and only if $T_{x}C$
is an isotropic (resp. coisotropic) subspace of $T_{x}M$ for every smooth
point $x\in C$.
\end{proposition}

Proof: For the isotropic case, the \textit{if} part of the assertion is
obvious. For the \textit{only if} part it suffices for us to check $T_{x}C$
being an isotropic subspace of $T_{x}M$ for a generic point $x\in C$ because
$\Omega|_{T_{x}C}=0$ is a closed condition among those $x$'s with constant
$\dim T_{x}C$. If $x$ is a smooth generic point in $C$ at which $\Omega
|_{T_{x}C}\neq0$, then
\[
\frac{\Omega\bar{\Omega}\omega^{k-2}}{\omega^{k}}|_{x}>0
\]
by the Riemann bilinear relation, where $k=\dim C$. Therefore we have
\[
\int_{C}\Omega\bar{\Omega}\omega^{k-2}>0
\]
violating $L_{\Omega}\left[  C\right]  =0$ assumption. Hence the claim. The
coisotropic case can be treated in a similar way. $\blacksquare$

\bigskip

Another corollary of the previous theorem is the following proposition which
generalizes Fujiki's observation that a generic complex structure in the
twistor family of $M$ has no curves or hypersurfaces.

\begin{proposition}
Suppose $M$ is a compact hyperk\"{a}hler manifold. For a general complex
structure in its twistor family, $M$ has no nontrivial isotropic or
coisotropic submanifold.
\end{proposition}

Proof: Suppose $C\subset M$ is an isotropic submanifold of positive dimension
$k$, then $\left[  C\right]  \cup\operatorname{Re}\Omega=\left[  C\right]
\cup\operatorname{Im}\Omega=0$ but $\left[  C\right]  \cup\omega^{k}>0$ in
$H^{\ast}\left(  M,\mathbb{R}\right)  $ because $C$ is a complex submanifold.
Other K\"{a}hler structures in the same twistor family can be written as
$a\operatorname{Re}\Omega+b\operatorname{Im}\Omega+c\omega$ for $\left(
a,b,c\right)  \in\mathbb{R}^{3}$ satisfying $a^{2}+b^{2}+c^{2}=1$. This
implies that $\left[  C\right]  $ can not be represented by an isotropic
complex submanifold in any other K\"{a}hler structures in this uncountable
family, except possibly $-\omega$ for $k$ even. On the other hand $\left[
C\right]  $ belongs to the integral cohomology of $M$, which is a countable
set. Hence we have our claim for the isotropic case. The coisotropic case can
also be argued in a similar way. $\blacksquare$

\bigskip

We have the following immediate corollary of the above proof.

\begin{corollary}
For any given $c\in H^{\ast}\left(  M,\mathbb{Z}\right)  $ in a compact
hyperk\"{a}hler manifold $M$, there is at most two complex structures in the
twistor family of $M$ such that $c$ can be represented by an isotropic or
coisotropic subvariety.
\end{corollary}

Note that isotropic or coisotropic submanifolds are plentiful when the whole
twistor family of complex structures on $M$ is considered. In the case of a K3
surface or an Abelian surface, they are complex curves and the numbers of such
form a beautiful generating function in terms of modular forms, as conjectured
by Yau and Zaslow in \cite{YZ} and proved by Bryan and the author in
\cite{BL1}, \cite{BL2} (also see \cite{BL3}).

\bigskip

\textbf{Coisotropic complete intersections}

When $M$ is projective, there are many ample hypersurfaces $C$ and they are
all coisotropic for trivial reasons. They are varieties of general type and
their Chern classes satisfy\footnote{These can be proven using the adjunction
formula, $c_{2k+1}\left(  M\right)  =0 $ and the Chern number inequality
$\int_{M}c_{2}\left(  M\right)  \left[  C\right]  ^{2n-2}\geq0$.} (i)
$c_{2k+1}\left(  C\right)  =c_{1}\left(  C\right)  c_{2k}\left(  C\right)  $
for all $k$ and (ii) $\int_{C}c_{2}\left(  C\right)  c_{1}\left(  C\right)
^{2n-3}\leq\int_{C}c_{1}\left(  C\right)  ^{2n-1}$.

The next theorem says that complete intersection coisotropic subvarieties can
be characterized by the vanishing of a Chern number. In particular
intersecting ample hypersurfaces will not give higher codimension coisotropic subvarieties.

\begin{theorem}
Suppose $M$ is a compact hyperk\"{a}hler manifold and $C$ is a complete
intersection subvariety of dimension $n+k$ with $k\leq n-2$. Then
\[
\int_{C}c_{2}\left(  M\right)  \left(  \Omega\bar{\Omega}\right)  ^{k}%
\omega^{n-k-2}\geq0\text{.}%
\]
Moreover the equality sign holds if and only if $C$ is coisotropic.
\end{theorem}

Proof: We can express the above quantity in term of the Bogomolov-Beauville
quadratic form $q$ by Fujiki's result (\cite{Fu}, see also the appendix) which
says for any $D_{1},\cdots,D_{2n-2}\in H^{2}\left(  M,\mathbb{C}\right)  $ we
have
\[
\int_{M}c_{2}\left(  M\right)  D_{1}\cdots D_{2n-2}=c\sum_{\substack{\left\{
i_{1},\cdots,i_{2n-2}\right\}  \\=\left\{  1,\cdots,2n-2\right\}  }}q\left(
D_{i_{1}},D_{i_{2}}\right)  \cdots q\left(  D_{i_{2n-3}},D_{i_{2n-2}}\right)
,
\]
for some constant $c$. When $D_{i}$'s are ample divisors, the left hand side
is strictly positive by the Chern number inequality for K\"{a}hler-Einstein
manifolds. On the other hand, $q\left(  D,D^{\prime}\right)  >0$ if $D$ is an
ample divisor and $D^{\prime}$ is an effective divisor. Therefore we have
$c>0$. More generally $q\left(  D,D^{\prime}\right)  \geq0$ when $D$ and
$D^{\prime}$ are effective divisors and they intersect transversely, moreover
the equality sign holds if and only if $\left(  \Omega\right)  ^{n-1}|_{D\cap
D^{\prime}}=0$. This is because
\[
q\left(  D,D^{\prime}\right)  =c^{\prime}\int_{M}\left[  D\right]  \left[
D^{\prime}\right]  \Omega^{n-1}\bar{\Omega}^{n-1}%
\]
for some explicit positive constant $c^{\prime}$.

We recall other basic properties of the Bogomolov-Beauville quadratic form:
$q\left(  \Omega,\bar{\Omega}\right)  >0$; $q\left(  \Omega,\omega\right)
=q\left(  \Omega,D\right)  =q\left(  \bar{\Omega},D\right)  =0$ for any
effective divisor $D$. Now $C$ is a complete intersection subvariety of $M$,
we write
\[
C=D_{1}\cap\cdots\cap D_{n-k}%
\]
for some effective divisors $D_{i}$'s. We have%

\[
\int_{C}c_{2}\left(  M\right)  \left(  \Omega\bar{\Omega}\right)  ^{k}%
\omega^{n-k-2}=\int_{M}c_{2}\left(  M\right)  \left[  D_{1}\right]
\cdots\left[  D_{n-k}\right]  \Omega^{k}\bar{\Omega}^{k}\omega^{n-k-2}.
\]
We apply Fujiki's result and above properties of $q$ and we get
\begin{align*}
&  \int_{C}c_{2}\left(  M\right)  \left(  \Omega\bar{\Omega}\right)
^{k}\omega^{n-k-2}\\
&  =cq\left(  \Omega,\bar{\Omega}\right)  ^{k}\sum q\left(  D_{i_{1}}%
,D_{i_{2}}\right)  q\left(  \omega,D_{i_{3}}\right)  \cdots q\left(
\omega,D_{i_{n-k}}\right) \\
&  +cq\left(  \omega,\omega\right)  q\left(  \Omega,\bar{\Omega}\right)
^{k}\sum q\left(  D_{i_{1}},D_{i_{2}}\right)  q\left(  D_{i_{3}},D_{i_{4}%
}\right)  q\left(  \omega,D_{i_{3}}\right)  \cdots q\left(  \omega,D_{i_{n-k}%
}\right)  +...\text{.}%
\end{align*}

Each term on the right hand side is non-negative. This implies
\[
\int_{C}c_{2}\left(  M\right)  \left(  \Omega\bar{\Omega}\right)  ^{k}%
\omega^{n-k-2}\geq0.
\]
Moreover it is zero if and only if $q\left(  D_{i},D_{j}\right)  =0$ for all
$i\neq j$. This is equivalent to $\Omega^{n-1}|_{D_{i}\cap D_{j}}=0$ for all
$i\neq j$ because $D_{i}\cap D_{j}$ is a complete intersection. We will prove
in the next lemma that this is equivalent to $\Omega^{k+1}|_{D_{1}\cap
\cdots\cap D_{n-k}}=0$, i.e. $C=D_{1}\cap\cdots\cap D_{n-k}$ is a coisotropic
subvariety of $M$. Hence the result. $\blacksquare$

\bigskip

The next lemma on linear algebra is needed in the proof of the above theorem
and it is also of independent interest.

\begin{lemma}
Let $M\cong\mathbb{C}^{2n}$ be a symplectic vector space with its symplectic
form $\Omega$ and $C$ is a codimension $m$ linear subspace in $M$. If we write
$C=D_{1}\cap\cdots\cap D_{m}$ for some hyperplanes $D_{i}$'s in $M$, then $C$
is coisotropic in $M$ if and only if $\Omega^{n-1}|_{D_{i}\cap D_{j}}=0$ for
all $i\neq j$.
\end{lemma}

Proof: We first prove the \textit{if }part by induction on the dimension of
$C$. The claim is trivial for $m=1$. When $m=2,$ it says that a codimension
two subspace $C$ in $M$ is coisotropic if $\Omega^{n-1}|_{C}=0$. This is true
and in general we have any codimension $m$ subspace $C$ in $M$ is coisotropic
if and only if $\Omega^{n-m+1}|_{C}=0$. By induction we assume $C=D_{1}%
\cap\cdots\cap D_{m+1}$ and the claim is true for $m$, i.e. $D_{1}\cap
\cdots\cap D_{m}$ is coisotropic. We can choose coordinates on $M$ such that
$\Omega=dz^{1}dz^{n+1}+\ldots+dz^{n}dz^{2n}$ and $D_{i}=\left\{
z^{i}=0\right\}  $ for $i=1,\ldots,m$. Suppose $%
{\textstyle\sum_{i=1}^{2n}}
a_{i}z^{i}=0$ is the defining equation for $D_{m+1}$. In order for $C$ to be
coisotropic, it suffices to show that $a_{i}=0$ for $n+1\leq i\leq n+m$. For
instance if $a_{n+1}\neq0$ then $\Omega^{n-1}|_{D_{1}\cap D_{m+1}}$ would be
nonzero. It is because $dz^{2}dz^{n+2}dz^{3}dz^{n+3}\cdots dz^{n}%
dz^{2n}|_{D_{1}\cap D_{m+1}}\neq0$ and all other summands of $\Omega^{n-1}$
would restrict to zero on $D_{1}\cap D_{m+1}$. If any other $a_{n+j}\neq0$ the
same argument applies by replacing $D_{1}$ with $D_{j}$. This contradicts our
assumption $\Omega^{n-1}|_{D_{i}\cap D_{j}}=0$ for all $i\neq j$. Therefore
$a_{i}=0$ when $n+1\leq i\leq n+m$ and hence the claim.

For the \textit{only if }part, we need to prove that $\Omega^{n-1}|_{D_{1}\cap
D_{2}}\neq0$ implies $\Omega^{n-l+1}|_{D_{1}\cap\cdots\cap D_{l}}\neq0$ for
any complete intersection $D_{1}\cap\cdots\cap D_{l} $. This is a simple
linear algebra exercise. Hence the lemma. $\blacksquare$

\bigskip

\subsection{\label{Sec Contract Birat}Contractions and birational transformations}

In this subsection we recall known results on the birational geometry of
hyperk\"{a}hler manifolds. The exceptional locus of any birational morphism is
always a coisotropic subvariety. The reduction and the projection of a
Lagrangian subvariety with respect to this coisotropic subvariety will be used
to obtain the Plucker type formula. Moreover the Legendre transformation that
we will discuss in section \ref{Sec Legendre} fits nicely into this picture.

In \textit{real }symplectic geometry, coisotropic submanifolds play a role in
the \textit{symplectic reductions} which produce symplectic manifolds of
smaller dimensions under favorable conditions (see for example section 5.1 of
\cite{BW}). For example when there is a Hamiltonian group action one can
construct a symplectic quotient by the reduction method.

In \textit{complex }symplectic geometry a coisotropic subvariety arises
naturally as the exceptional set of any birational contraction, and its
geometric structures can even be described rather explicitly. Roughly speaking
it is a generic $\mathbb{P}^{k}$-bundle over another symplectic manifold of
dimension $2\left(  n-k\right)  $.\footnote{When $k=1$ the generic fiber can
also be an ADE configuration of $\mathbb{P}^{1}$.} When the exceptional locus
is a honest $\mathbb{P}^{k}$-bundle, Mukai \cite{Mu1} showed that one can
produce another symplectic manifold of the \textit{same }dimension by
replacing the coisotropic $\mathbb{P}^{k}$-bundle by its dual $\mathbb{P}^{k}%
$-bundle, the Mukai elementary modification. In dimension four this turns out
to be the only possible way to generate birational transformations provided
certain normality conditions hold (\cite{BHL}).

All these structures will give important functors among categories of
Lagrangian subvarieties that we will define in later sections.

\bigskip

\textbf{Characteristic foliations and reductions}

For any submanifold $C$ in $M$, there is an exact sequence relating their
tangent bundles. Using the isomorphism between $T_{M}$ and $T_{M}^{\ast}$
induced from the symplectic form $\Omega$, we have the following exact
diagram,
\[%
\begin{array}
[c]{ccccccc}%
0\rightarrow &  T_{C} & \rightarrow &  T_{M}|_{C} & \rightarrow &  N_{C/M} &
\rightarrow0\\
& \downarrow &  & \,\updownarrow\cong &  & \uparrow & \\
0\leftarrow &  T_{C}^{\ast} & \leftarrow &  T_{M}^{\ast}|_{C} & \leftarrow &
N_{C/M}^{\ast} & \leftarrow0,
\end{array}
\]
where $N_{C/M}$ is the normal bundle of $C$ inside $M$. The coisotropic
condition on $C$ is equivalent to the triviality of the composition
homomorphism $N_{C/M}^{\ast}\rightarrow N_{C/M}$. When this happens we have an
injective homomorphism $N_{C/M}^{\ast}\rightarrow T_{C}$ and we denote the
quotient bundle as $S$,
\[
0\rightarrow N_{C/M}^{\ast}\rightarrow T_{C}\rightarrow S\rightarrow0\text{.}%
\]

A different viewpoint is that the restriction of $\Omega$ to $T_{C}$ is no
longer non-degenerate, and the kernel is precisely given by $N_{C/M}^{\ast}$.
It is then clear that the bundle $S$ inherits a natural symplectic structure.
In particular $c_{2k+1}\left(  S\right)  =0$ for all $k$. This implies that
$T_{C}$ cannot be a stable bundle with nonnegative degree with respect to any
polarization. For instance a coisotropic submanifold in $M$ can never be a
Calabi-Yau manifold or an irreducible hyperk\"{a}hler manifold unless it is a
Lagrangian. If $C$ is an Abelian variety, then the above exact sequence splits
and both $N_{C/M}^{\ast}$ and $S$ are trivial bundles over $C$. This happens
in Lagrangian fibrations on $M$.

\bigskip

The distribution on $C$ given by $N_{C/M}^{\ast}\subset T_{C}$ is always
integrable (see e.g. p.67 of \cite{BW}). This is called the
\textit{characteristic foliation }of the coisotropic submanifold. When the
\textit{leaf space} $B$ is a smooth manifold, it carries a natural symplectic
form, moreover this form pulls back to the symplectic structure on the bundle
$S$ over $C$. This symplectic manifold $B$ if exists, is called the
\textit{reduction} of the coisotropic submanifold $C$ (see section \ref{Sec
Red Proj functor} for further details).

Roughly speaking what we did is to contract degenerate directions in a
coisotropic submanifold to construct the symplectic form on the quotient space
$B$. A natural question is whether we can perform this inside $M$ rather than
$C$. This would correspond to having a birational contraction of $M$.

\bigskip

\textbf{Birational contraction }

Suppose that
\[
\pi:M\rightarrow Z
\]
is a projective birational morphism, which we call a\textit{\ birational
contraction}. Let $C$ be the exceptional locus in $M$ and $B=\pi\left(
C\right)  $. Using the Grauert-Riemenschneider vanishing result, one can show
that $Z$ has at worst rational singularities and $\pi$ is a crepant resolution
with fibers uniruled varieties. In particular every fiber $F$ of $\pi
_{C}:C\rightarrow B$ has no nontrivial holomorphic form and therefore $F$ is
isotropic inside $M$. Similarly $\left(  \iota_{v}\Omega\right)  |_{F}=0$ for
any $v\in TB$ and this implies that $T_{vert}C\subseteq\left(  TC\right)
^{\bot}.$

Wierzba, Namikawa and Hu-Yau (\cite{Wi}, \cite{Na}, \cite{HY}) show that
generically,
\[
T_{vert}C=\left(  TC\right)  ^{\bot}.
\]
This implies that $C$ is a coisotropic subvariety and its generic
characteristic foliation is given by fibers of $\pi_{C}$, in particular $B$
has a natural symplectic structure outside its singular locus. The proof goes
roughly as follow: It suffices to show that $\dim F\geq\dim M-\dim C$; we
consider deformations of rational curves in $C$ which cannot bend and break.
On the one hand any such curve cannot move away from the contractible locus
$C$ and therefore the dimension of its deformation space cannot be too big
compare to $\dim F$ by a result of Cho and Miyaoka. On the other hand, by a
result of Ran, every rational curve in $M$ has at least $2n-2$ deformation
directions even though the expected dimension is only $2n-3$. The fact that
generic complex structures in the twistor family of $M$ has no curve is
responsible for the extra dimension. Combining this we have $\dim F\geq\dim
M-\dim C$ and therefore $T_{vert}C=\left(  TC\right)  ^{\bot}.$

In fact each smooth\footnote{Normality is already sufficient by the work of
Shepherd-Barron \cite{Sh}.} irreducible component of $F$ is a projective space
$\mathbb{P}^{k}$. The proof of this combines above arguments with
Cho-Miyaoka's characterization of $\mathbb{P}^{k}$ as the only smooth variety
with the property that every rational curve moves in at least $2k-2$
dimensional family.

Note that $\pi_{C}:C\rightarrow B$ is not necessary a honest bundle, the
dimension of fibers can jump and $B$ could be singular too. There are some
structure theorems describing such behaviors and we expect them to admit good
stratifications as defined by Markman \cite{Mar}.

\bigskip

\textbf{Mukai elementary modification}

Using previous arguments, one can show that every $\mathbb{P}^{k}$-bundle
\[
\pi_{C}:C\rightarrow B
\]
inside $M$ has dimension less than or equal to $2n-k$. When $\dim C=2n-k$, $C
$ is coisotropic and its characteristic foliation is precisely given by fibers
of $\pi_{C}$. When this happens Mukai \cite{Mu1} constructs \textit{another}
symplectic manifold $M^{\prime}$ by replacing each fiber $\mathbb{P}^{k}$ with
its dual projective space $\left(  \mathbb{P}^{k}\right)  ^{\ast}$. This is
called the \textit{Mukai elementary modification.}

When $C=\mathbb{P}^{n}$ it has a neighborhood which is symplectomorphic to a
neighborhood of the zero section in $T^{\ast}\mathbb{P}^{n}$. In this case
this modification is simply the birational map $\Phi$ in the introduction. We
will study a Legendre transformation along any $\mathbb{P}^{n}$ in $M$ in the
section \ref{Sec Legendre}. The general case can be regarded as a family
version of this.

\section{Lagrangian submanifolds}

Lagrangian submanifolds are the smallest coisotropic submanifolds or the
biggest isotropic submanifolds. They are the most important objects in
symplectic geometry, 'Everything is a Lagrangian manifold' as Weinstein
described their roles in real symplectic geometry.

In this section we will first establish some basic topological and geometrical
properties of Lagrangian submanifolds in a hyperk\"{a}hler manifold. We will
explain the relationship between the second fundamental form and the special
K\"{a}hler structure on the moduli space of Lagrangian submanifolds. Then we
will study a Lagrangian category and intersection theory of Lagrangian
subvarieties, a reduction functor and a projection functor. These structures
will be needed to study the Legendre transformation and the Plucker type
formula in the next section.

\subsection{Properties of Lagrangian submanifolds}

When $C$ is a Lagrangian submanifold of $M$, then $N_{C/M}\cong T_{C}^{\ast} $
and we have the following exact sequence
\[
0\rightarrow T_{C}\rightarrow T_{M}|_{C}\rightarrow T_{C}^{\ast}%
\rightarrow0\text{.}%
\]
By the Whitney sum formula we have
\begin{align*}
\iota^{\ast}c\left(  T_{M}\right)   &  =c\left(  T_{C}\oplus T_{C}^{\ast
}\right) \\
&  =c\left(  T_{C}\otimes_{\mathbb{R}}\mathbb{C}\right)  \text{,}%
\end{align*}
where $\iota:C\rightarrow M$ is the inclusion morphism.

This implies that Pontrjagin classes of $C$ are determined by Chern classes of
$M$ as follows,
\[
p_{k}\left(  C\right)  =\left(  -1\right)  ^{k}\iota^{\ast}c_{2k}\left(
M\right)  .
\]

In particular, Pontrjagin numbers of $C$ depends only on the cohomology class
$\left[  C\right]  \in H^{2n}\left(  M,\mathbb{Z}\right)  $. By the celebrated
theorem of Thom, this determines the rational cobordism type of $C $ and we
have proven the following theorem.

\begin{theorem}
Given a fixed cohomology class $c\in H^{2n}\left(  M,\mathbb{Z}\right)  $, it
determines the rational cobordism class of any possible Lagrangian submanifold
in $M$ representing $c$.
\end{theorem}

In particular, the signature of the intersection product on $H^{\ast}\left(
C,\mathbb{R}\right)  $ is also determined by $c$. Using the Hirzebruch
signature formula we have
\[
Signature\left(  C\right)  =\int_{M}c\cup\sqrt{L_{M}}\text{.}%
\]
Using $N_{C/M}\cong T_{C}^{\ast}$ we also have the following formula for the
Euler characteristic of $C$,%

\[
\chi\left(  C\right)  =\left(  -1\right)  ^{n}\int_{M}c\cup c.
\]
Similarly if $C$ is spin, the \^{A}-genus is given by
\[
\hat{A}\left[  C\right]  =\int_{M}c\cup\sqrt{Td_{M}}.
\]
The reason is when restricting to $C$, using the multiplicative property of
the $\hat{A}$-genus, we have
\[
Td_{M}=\hat{A}\left(  T_{M}\right)  =\hat{A}\left(  T_{C}+T_{C}^{\ast}\right)
=\hat{A}\left(  T_{C}\right)  \hat{A}\left(  T_{C}^{\ast}\right)  =\hat
{A}\left(  T_{C}\right)  ^{2}.
\]

As an immediate corollary of this and the standard Bochner arguments, we can
show that if $C$ is an even dimensional spin Lagrangian submanifold of $M$,
then $T_{C}$ cannot have positive scalar curvature unless $\int_{C}%
\sqrt{Td_{M}}=0.$

\bigskip

Recall from section \ref{Sec Iso-Coiso property} that we have following
characterizations of a $n$ dimensional submanifold $C$ in $M$ being a
Lagrangian. First $C$ is Lagrangian if $H^{2,0}\left(  C\right)  =0$. Second
if $C$ is a complete intersection then
\[
\int_{C}c_{2}\left(  M\right)  \omega^{n-2}\geq0\text{.}%
\]
Moreover the equality sign holds if and only if $C$ is a Lagrangian. When this
happens the tangent bundle of $C$ splits into direct sum of line bundles, this
is a strong constraint upon $C$. For example when $n=2$ the signature of $C$
would have to be zero. Third $C$ is Lagrangian if and only if its Poincar\'{e}
dual is a $\Omega$-primitive class, i.e. $\Lambda_{\Omega}\left[  C\right]
=0$, or equivalently, $L_{\Omega}\left[  C\right]  =0$\textit{. }This implies
that Lagrangian property is invariant under deformations of $C$ and it also
allows us to define singular, and possibly non-reduced, Lagrangian
subvarieties of $M$.

\bigskip

Remark: When $C$ is a Lagrangian submanifold, then $\left[  C\right]  \in
H^{n,n}\left(  M,\mathbb{Z}\right)  \cap Ker\left(  L_{\Omega}\right)  $. In
fact being a class of type $\left(  n,n\right)  $ follows from the $\Omega
$-primitivity. More precisely, we have
\begin{align*}
Ker\left(  L_{\Omega}\right)  \cap\Omega^{2n}\left(  M,\mathbb{R}\right)   &
\subset\Omega^{n,n}\left(  M\right)  ,\\
Ker\left(  L_{\Omega}\right)  \cap H^{2n}\left(  M,\mathbb{R}\right)   &
\subset H^{n,n}\left(  M\right)  .
\end{align*}
The proof of these inclusions simply uses the hard Lefschetz decomposition of
$\Omega^{\ast,\ast}$ for the $\mathbf{sl}\left(  2\right)  $\textbf{-}action
generated by $L_{\Omega}$ and $\Lambda_{\Omega}$. A similar result of Hitchin
\cite{Hi2} says that if $C$ is a \textit{real} submanifold of $M$ which is a
\textit{real} Lagrangian submanifold with respect to both $\operatorname{Re}%
\Omega$ and $\operatorname{Im}\Omega$, then $C$ is a \textit{complex
}submanifold and therefore a Lagrangian submanifold of $M$ with respect to
$\Omega$.

\bigskip

\textbf{Some examples}

(1) When $M$ is a K3 surface or an Abelian surface, any curve $C$ in $M$ is
Lagrangian for dimensional reason. The Hilbert scheme of $n$ points in $M$,
denote $S^{\left[  n\right]  }M$, is again a hyperk\"{a}hler manifold as shown
by Fujiki and Beauville. Moreover $S^{\left[  n\right]  }C\subset S^{\left[
n\right]  }M$ is a Lagrangian submanifold. Notice that $S^{\left[  n\right]
}C$ is simply the symmetric product of the curve $C$, i.e. $S^{n}C$.

\bigskip

(2) For any K\"{a}hler manifold $X$, its cotangent bundle $T^{\ast}X$ carries
a natural symplectic structure. Moreover the conormal bundle of any
submanifold of $X$ is a Lagrangian submanifold of $T^{\ast}X$. In fact they
can be characterized as those closed submanifolds of $T^{\ast}X$ which are
invariant under the natural $\mathbb{C}^{\times}$-action of scaling along
fibers on $T^{\ast}X$. They are called \textit{conical submanifolds.}

\bigskip

\subsection{Second fundamental form and moduli space}

In this subsection we study the differential geometry of a Lagrangian
submanifold $C$ in $M$. The second fundamental form of a submanifold is a
symmetric two tensor on $C$ with valued in the normal bundle and it is defined
as follow, $X\otimes Y\rightarrow\left(  \nabla_{X}Y\right)  ^{N}$ for any
tangent vectors $X,Y$ on $C$. Here $\nabla$ is the Levi-Civita connection on
$M$ and $\left(  \bullet\right)  ^{N}$ is the orthogonal projection to normal directions.

The Lagrangian condition gives us a symmetric three tensor and it defines a
cubic form on $H^{0}\left(  C,T_{C}^{\ast}\right)  $. By varying $C$ this
cubic form will determine a special K\"{a}hler structure on the moduli space
of Lagrangian submanifolds. We will also give a brief physical explanations of
such structure using supersymmetry.

\subsubsection{Second fundamental form and a cubic form}

If $C$ is a Lagrangian submanifold of $M$, then $T_{M}|_{C}$ is a symplectic
vector bundle with a Lagrangian subbundle $T_{C}$.

\bigskip

\textbf{Extension class of a Lagrangian subbundle}

Suppose $E$ is a symplectic vector bundle over a complex manifold $C$. Given
any Lagrangian subbundle $S$, by contracting the symplectic form $\Omega$ with
elements in $S$, we have a naturally identification between quotient bundle
$E/S$ and the dual bundle $S^{\ast}$ of $S$. Therefore we have the following
short exact sequence,
\[
0\rightarrow S\rightarrow E\rightarrow S^{\ast}\rightarrow0\text{.}%
\]
We consider the image of $1_{S^{\ast}}$ under the homomorphism $Hom_{O_{C}%
}\left(  S^{\ast},S^{\ast}\right)  \rightarrow Ext_{O_{C}}^{1}\left(  S^{\ast
},S\right)  $, coming from the corresponding long exact sequence, and call it
$\alpha$. This is the extension class
\[
\alpha\in Ext_{O_{C}}^{1}\left(  S^{\ast},S\right)  \cong H^{1}\left(
C,S\otimes S\right)  ,
\]
associated to the above short exact sequence. The next lemma says $\alpha$
lies inside the symmetric tensor product of $S$.

\begin{lemma}
If $S$ is a Lagrangian subbundle of a symplectic vector bundle $E$ over $C$,
then the corresponding extension class $\alpha$ lies inside $H^{1}\left(
C,Sym^{2}S\right)  $.
\end{lemma}

Proof: For any \textit{non-holomorphic} splitting $E=S\oplus S^{\ast}$ of the
above exact sequence, we write the $\bar{\partial}$-operator of $E\;$in terms
of this decomposition,
\[
\bar{\partial}_{E}=\left(
\begin{array}
[c]{cc}%
\bar{\partial}_{S} & A\\
0 & \bar{\partial}_{S^{\ast}}%
\end{array}
\right)
\]
then the off-diagonal term $A\in\Omega^{0,1}\left(  C,S\otimes S\right)  $
represents the extension class $\alpha$. On $S\oplus S^{\ast}$, there is a
canonical symplectic form constructed by the natural pairing between $S$ and
$S^{\ast}$. Using linear algebra we can choose a splitting suitably such that
$\Omega$ on $E$ becomes the canonical symplectic form on $S\oplus S^{\ast}$.
It is easy to verify that $\bar{\partial}_{E}\Omega=0$ is then equivalent to
$A\in\Omega^{0,1}\left(  C,Sym^{2}S\right)  $. Therefore we have $\alpha\in
H^{1}\left(  C,Sym^{2}S\right)  $. $\blacksquare$

\bigskip

Given any Hermitian metric on $E$ there is unique Hermitian connection $D_{E}
$ on $E$ satisfying $\left(  D_{E}\right)  ^{0,1}=\bar{\partial}_{E}$. Using
the metric we can also decompose $E$ into $S\oplus S^{\ast}$ as an orthogonal
decomposition. The above differential form $A\in\Omega^{0,1}\left(  C,S\otimes
S\right)  $ is simply the second fundamental form of the subbundle $S\subset
E$ in the context of differential geometry. If the symplectic form $\Omega$ on
$E$ becomes the canonical symplectic form on the orthogonal decomposition
$S\oplus S^{\ast}$ then the above proof gives $A\in\Omega^{0,1}\left(
C,Sym^{2}S\right)  $. To achieve this we need the Hermitian metric $h$ on $E$
to be compatible with $\Omega$ in the sense that if we define operators $I$
and $K$ by the following,
\[
\Omega\left(  v,w\right)  =h\left(  Iv,w\right)  +ih\left(  Kv,w\right)
\text{,}%
\]
then $I,J,K$ defines a fiberwise quaternionic structure on $E$, i.e.
$I^{2}=J^{2}=K^{2}=IJK=-1$.

\begin{proposition}
Suppose that $E$ is a symplectic bundle over $C$ with compatible quaternionic
structure. Then the second fundamental form of any Lagrangian subbundle
$S\subset E$ lies inside $\Omega^{0,1}\left(  C,Sym^{2}S\right)  $.
\end{proposition}

Proof: Let us denote the orthogonal complement of $S$ in $E$ as $T$. One can
verify using the quaternionic structure on $E$ that $T=I\left(  S\right)  $.
This implies that $T$ is a Lagrangian subbundle of $E$. Therefore the
symplectic form $\Omega$ on $E$ can be identified as the canonical symplectic
form on $S\oplus S^{\ast}$ via the natural identification $T\cong E/S\cong
S^{\ast}$. By the proof of the previous lemma, we have our proposition.
$\blacksquare$

\bigskip

\textbf{A cubic form}

When $C$ is a Lagrangian submanifold in $M$, then its tangent bundle $S=T_{C}
$ is a Lagrangian subbundle of the symplectic bundle $E=T_{M}|_{C}$. The
hyperk\"{a}hler metric on $M$ does give us a compatible quaternionic structure
on $T_{M}|_{C}$. Therefore the second fundamental form $A$ satisfies,
\[
A\in\Omega^{0,1}\left(  C,Sym^{2}T_{C}\right)  \text{.}%
\]

On the other hand the induced metric on $C$ determines a natural
identification between $\left(  T_{C}^{\ast}\right)  ^{0,1}$ and $\left(
T_{C}\right)  ^{1,0}=T_{C}$. Hence we can identify $A$ as an element in
$\Gamma\left(  C,T_{C}\otimes Sym^{2}T_{C}\right)  $ and we denote it as
$\widetilde{A}$. Explicitly if $z^{i}$'s is a local coordinate system on $C$
and the induced metric has components $g_{i\bar{j}}$'s, then
\begin{align*}
\widetilde{A} &  =%
{\textstyle\sum}
\widetilde{A}^{ijk}\frac{\partial}{\partial z^{i}}\otimes\frac{\partial
}{\partial z^{j}}\otimes\frac{\partial}{\partial z^{k}},\\
\widetilde{A}^{ijk} &  =%
{\textstyle\sum}
g^{i\bar{l}}\left(  A^{jk}\right)  _{\bar{l}}\text{.}%
\end{align*}
From basic differential geometry we know that the second fundamental form of
any submanifold is a symmetric tensor with valued in the normal bundle, this
reflects the fact the Levi-Civita connection on $M$ is torsion free. As a
consequence $\widetilde{A}^{ijk}$ is completely symmetric in $i,j$ and $k$,
i.e. $\widetilde{A}$ is a trilinear symmetric multi-vector field on $C$,%

\[
\widetilde{A}\in\Gamma\left(  C,Sym^{3}T_{C}\right)  \text{.}%
\]

We use $\widetilde{A}$ to define a natural cubic form on $H^{0}\left(
C,T_{C}^{\ast}\right)  $ as follows.

\begin{definition}
If $M$ is a hyperk\"{a}hler manifold and $C$ is a compact Lagrangian
submanifold of $M$ then we define a cubic form on $H^{0}\left(  C,T_{C}^{\ast
}\right)  $,
\begin{align*}
\mathbf{c}_{C}:Sym^{3}H^{0}\left(  C,T_{C}^{\ast}\right)  \rightarrow
\mathbb{C}\\
\mathbf{c}_{C}\left(  \phi,\eta,\zeta\right)  =\int_{C}\phi_{i}\eta_{j}%
\zeta_{k}\widetilde{A}^{ijk}\frac{\omega^{n}}{n!}\text{.}%
\end{align*}
\end{definition}

Remark: The above cubic form $\mathbf{c}_{C}$ coincides with the composite of
natural homomorphisms
\[
H^{0}\left(  C,T_{C}^{\ast}\right)  ^{\otimes3}\overset{\alpha}{\rightarrow
}H^{1}\left(  C,O_{C}\right)  \otimes H^{0}\left(  C,T_{C}^{\ast}\right)
\overset{\otimes}{\rightarrow}H^{1}\left(  C,T_{C}^{\ast}\right)
=H^{1,1}\left(  C\right)  \overset{\Lambda}{\rightarrow}\mathbb{C}\text{.}%
\]
In particular it depends only on the extension class $\alpha$.

\subsubsection{Moduli space of Lagrangian submanifolds}

If $C$ is Lagrangian submanifold in $M$ then any deformation of $C$ inside $M
$ remains a Lagrangian. Therefore the moduli space of Lagrangian submanifolds
is the same as the moduli space of submanifolds representing the same
cohomology class in $M$. We denote it by $\mathcal{M}$. This moduli space is
always smooth, namely any infinitesimal deformation of $C$ inside $M $ is
always unobstructed. This can be proven by applying a twistor rotation on
McLean's result \cite{Mc} on unobstructedness of deformations of special
Lagrangian submanifolds in Calabi-Yau manifolds. This can also be proved by
algebraic geometric means using T$^{1}$-lifting method developed by Ran.

We recall that infinitesimal deformations of $C$ in $M$ are parametrized by
$H^{0}\left(  C,T_{C}^{\ast}\right)  $ because the Lagrangian condition gives
a natural identification between $N_{C/M}$ and $T_{C}^{\ast}$. It is not
difficult to see that contracting with the extension class $\alpha\in
H^{1}\left(  C,Sym^{2}T_{C}\right)  $ gives a homomorphism
\[
\iota_{\alpha}:H^{0}\left(  C,T_{C}^{\ast}\right)  \rightarrow H^{1}\left(
C,T_{C}\right)
\]
which associates to an infinitesimal deformation of $C$ in $M$ to the
corresponding infinitesimal deformation of the complex structure on $C$. Next
we will show that this cubic form $\mathbf{c}$, constructed from the second
fundamental form, determines the special K\"{a}hler structure on $\mathcal{M}$.

\bigskip

\textbf{Moduli space as a special K\"{a}hler manifold}

We first give a physical reason for the existence of a natural special
K\"{a}hler structure on $\mathcal{M}$: $\sigma$-model studies maps from
Riemann surfaces to a fix target manifold. When the target manifold $M$ is
hyperk\"{a}hler, it is well-known that its $\sigma$-model has $N=4$
supersymmetry (or SUSY). If domain Riemann surfaces have boundary components
then we would require their images lie inside a fix Lagrangian submanifold
$C\subset M$. In this case only half of the SUSY can be preserved and we have
a $N=2$ SUSY theory. Moduli space of such theories are generally known
(physically) to process special K\"{a}hler geometry. In our situation, this is
simply the moduli space of Lagrangian submanifolds in $M$.

\bigskip

Let us recall the definition of a special K\"{a}hler manifold (see \cite{Fr}
for details). A special K\"{a}hler structure on a K\"{a}hler manifold
$\mathcal{M}$ is carried by a holomorphic cubic tensor
\[
\Xi\in H^{0}\left(  \mathcal{M},Sym^{3}T_{\mathcal{M}}^{\ast}\right)  .
\]
Using the K\"{a}hler metric $g_{\mathcal{M}}$ on $\mathcal{M}$ we can identify
$\Xi$ with a tensor $\mathbf{A}\in\Omega^{1,0}\left(  \mathcal{M}%
,EndT_{\mathbb{C}}\mathcal{M}\right)  $ as follows
\[
\Xi=-\omega_{\mathcal{M}}\left(  \pi^{1,0},\left[  \mathbf{A},\pi
^{1,0}\right]  \right)
\]
where $\pi^{1,0}\in\Omega^{1,0}\left(  T_{\mathbb{C}}\mathcal{M}\right)  $ is
constructed from the inclusion homomorphism $T^{1,0}\subset T_{\mathbb{C}} $.
In terms of local coordinates we have $\left(  \mathbf{A}_{l}\right)
_{j}^{\bar{k}}=i\Xi_{ljm}g_{\mathcal{M}}^{m\bar{k}}$. If we denote the
Levi-Civita connection on $\mathcal{M}$ as $\nabla^{LC}$, then the special
K\"{a}hler condition is $\nabla=\nabla^{LC}+\mathbf{A}+\mathbf{\bar{A}}$
defines a torsion free flat symplectic connection on the tangent bundle and it
satisfies $\nabla\wedge J=0$.

\bigskip

When $\mathcal{M}$ is the moduli space of Lagrangian submanifolds in $M$,
Hitchin \cite{Hi} shows that it has a natural special K\"{a}hler structure. We
are going to show that this special K\"{a}hler structure on $\mathcal{M}$ is
given by the previous cubic form.

At any given point $\left[  C\right]  \in\mathcal{M}$ with $C$ a Lagrangian
submanifold in $M$, we have a natural identification of the tangent space
$T_{\mathcal{M},\left[  C\right]  }=H^{0}\left(  C,T_{C}^{\ast}\right)  $.
From previous discussions we have a cubic form on $H^{0}\left(  C,T_{C}^{\ast
}\right)  $.
\[
\mathbf{c}_{C}:Sym^{3}H^{0}\left(  C,T_{C}^{\ast}\right)  \rightarrow
\mathbb{C}\text{.}%
\]
By varying the point $\left[  C\right]  $ in $\mathcal{M}$ this determines a
cubic tensor,
\[
\mathbf{c}\in\Gamma\left(  \mathcal{M},Sym^{3}T_{\mathcal{M}}^{\ast}\right)  .
\]

\begin{theorem}
Suppose $\mathcal{M}$ is the moduli space of Lagrangian submanifolds in a
compact hyperk\"{a}hler manifold $M$. Then $\mathbf{c}\in H^{0}\left(
\mathcal{M},Sym^{3}T_{\mathcal{M}}^{\ast}\right)  $ and it determines a
special K\"{a}hler structure on $\mathcal{M}$.
\end{theorem}

Proof: We are going to prove this theorem by identifying the cubic form in the
special K\"{a}hler structure on $\mathcal{M}$ with the above cubic form
$\mathbf{c}$.

The cubic form that determines the special K\"{a}hler structure on
$\mathcal{M}$ was described by Donagi and Markman in \cite{DM}. First
$\mathcal{M}$ can be viewed it as the base space of the moduli space of
universal compactified Jacobian $\mathcal{J}$ for Lagrangian submanifolds in
$M$. The space $\mathcal{J}$ has a natural holomorphic symplectic structure
and the natural morphism $\mathcal{J}$ $\rightarrow\mathcal{M}$ is a
Lagrangian fibration. We consider the variation of Hodge structures, i.e.
periods, for this family of Abelian varieties over $\mathcal{M}$. Its
differential at a point $\left[  C\right]  \in\mathcal{M}$ gives a
homomorphism,
\[
T_{\mathcal{M},\left[  C\right]  }\rightarrow S^{2}V,
\]
where
\begin{align*}
V  & =H^{1}\left(  Jac\left(  C\right)  ,O_{Jac\left(  C\right)  }\right)
=H^{1}\left(  C,O_{C}\right)  ,\\
T_{\mathcal{M},\left[  C\right]  }  & =H^{0}\left(  C,T_{C}^{\ast}\right)  .
\end{align*}
Using the induced K\"{a}hler structure on $C$, we can identify $H^{1}\left(
C,O_{C}\right)  ^{\ast}$ with $H^{0}\left(  C,T_{C}^{\ast}\right)  $. We thus
obtain a homomorphism
\[%
{\textstyle\bigotimes^{3}}
T_{\mathcal{M},\left[  C\right]  }\rightarrow\mathbb{C}%
\]
This tensor is a symmetric tensor. The corresponding cubic form $\Gamma\left(
\mathcal{M},Sym^{3}T_{\mathcal{M}}^{\ast}\right)  $ is the one that determines
the special K\"{a}hler structure on $\mathcal{M}$ (see \cite{Fr}).

The above homomorphism $T_{\mathcal{M},\left[  C\right]  }\rightarrow S^{2}V$
can be identified as the composition of the natural homomorphism $H^{0}\left(
C,N_{C/M}\right)  \rightarrow H^{1}\left(  C,T_{C}\right)  $ and the natural
homomorphism between variation of complex structures on $C$ and on its
Jacobian variety,
\[
H^{1}\left(  C,T_{C}\right)  \rightarrow H^{1}\left(  Jac\left(  C\right)
,T_{Jac\left(  C\right)  }\right)  \text{.}%
\]
From standard Hodge theory, an Abelian variety is determined by its period, or
equivalently its weight one Hodge structure, and its variation is simply given
by the usual cup product homomorphism,
\[
H^{1}\left(  C,T_{C}\right)  \otimes H^{1,0}\left(  C\right)  \overset{\cup
}{\rightarrow}H^{0,1}\left(  C\right)  ,
\]
or equivalently,
\begin{gather*}
H^{1}\left(  C,T_{C}\right)  \otimes H^{0}\left(  C,T_{C}^{\ast}\right)
\overset{\cup}{\rightarrow}H^{1}\left(  C,O_{C}\right) \\
\left(  \phi_{\bar{j}}^{i}\frac{\partial}{\partial z^{i}}\otimes d\bar{z}%
^{j}\right)  \otimes\left(  \beta_{k}dz^{k}\right)  \rightarrow\phi_{\bar{j}%
}^{i}\beta_{i}d\bar{z}^{j}.
\end{gather*}
Recall under the identification $H^{0}\left(  C,N_{C/M}\right)  \cong
H^{0}\left(  C,T_{C}^{\ast}\right)  $ we have
\begin{gather*}
H^{0}\left(  C,T_{C}^{\ast}\right)  \rightarrow H^{1}\left(  C,T_{C}\right) \\
\alpha_{i}dz^{i}\rightarrow\alpha_{i}A_{\bar{j}}^{ik}\frac{\partial}{\partial
z^{k}}\otimes d\bar{z}^{j},
\end{gather*}
where $A=\Sigma A_{\bar{j}}^{ik}\frac{\partial}{\partial z^{i}}\otimes
\frac{\partial}{\partial z^{k}}\otimes d\bar{z}^{j}$ is the second fundamental
form of $C$ in $M$. Therefore the corresponding cubic form $\otimes^{3}%
H^{0}\left(  C,T_{C}^{\ast}\right)  \rightarrow\mathbb{C}$ is given by
\[
\left(  \alpha,\beta,\gamma\right)  \rightarrow\int_{C}\alpha_{i}A_{\bar{j}%
}^{ik}\beta_{k}\gamma_{l}g^{l\bar{j}}\frac{\omega^{n}}{n!}\text{.}%
\]
This is precisely the cubic form $\mathbf{c}_{C}$. Hence the result.
$\blacksquare$

\bigskip

Remark: On a special K\"{a}hler manifold $\mathcal{M}$ its Riemannian
curvature tensor is given by
\[
R_{i\bar{j}k\bar{l}}=-g_{\mathcal{M}}^{p\bar{q}}\Xi_{ijp}\bar{\Xi}%
_{jlq}\text{.}%
\]
This implies that the Ricci curvature of $\mathcal{M}$ is non-negative and the
scalar curvature equals $4\left|  \Xi\right|  ^{2}\geq0$. Combining this with
our earlier discussions on the second fundamental form, we show that the
special K\"{a}hler metric on the moduli space of Lagrangian submanifolds is
\textit{flat} at $\left[  C\right]  \in\mathcal{M}$ if the natural exact
sequence $0\rightarrow T_{C}\rightarrow T_{M}|_{C}\rightarrow N_{C/M}%
\rightarrow0$ splits.

\subsection{Lagrangian category}

On a \textit{real} symplectic manifold $\left(  M,\omega\right)  $, Fukaya
proposed a category of Lagrangian submanifolds. The space of morphisms between
two Lagrangian submanifolds $L_{1},L_{2}$ is the Floer cohomology group
$HF\left(  L_{1},L_{2}\right)  $. It is defined in terms of the number of
holomorphic disks (i.e. instantons) bounding $L_{1}$ and $L_{2}$. The
dimension of the space of such holomorphic disks can be computed using the
index theorem, and expressed in terms of their Maslov indexes. Kontsevich
conjectured that Fukaya category on a Calabi-Yau manifold is equivalent to the
derived category of coherent sheaves of the \textit{mirror }Calabi-Yau
manifold, the so-called homological mirror symmetry conjecture. In this paper
we study the holomorphic analog of the Fukaya category. The definition of this
category is probably well-known to experts in this subject. We will establish
some basic properties of it and they will play an important role in the
discussions of the Pl\"{u}cker type formula in section \ref{Sec Legendre}.

\subsubsection{Definition of the category}

\textbf{Non-existence of instanton corrections}

Before we give the definition of the Lagrangian category we first demonstrate
the absent of \textit{instanton corrections }in the hyperk\"{a}hler setting.
Suppose $M$ is a hyperk\"{a}hler manifold with a preferred complex structure
$J$ among $I,J,K$. When $C$ is a Lagrangian submanifold in $M$ with respect to
the symplectic structure $\Omega=\omega_{I}+i\omega_{K}$ then it is a
\textit{real} Lagrangian submanifold of the \textit{real} symplectic manifold
$\left(  M,\omega_{\theta}\right)  $ with $\omega_{\theta}=\cos\theta
\omega_{I}+\sin\theta\omega_{K}$ for any $\theta\in\lbrack0,2\pi)$. For any
fixed $\theta$,\textit{\ instantons }means $J_{\theta}$-holomorphic disks
bounding $C$, where $J_{\theta}=\cos\theta I+\sin\theta K$. In Physics they
contribute to the \textit{correlation functions}, which is independent of
$\theta$ in any TQFT. In particular there is no instanton effects if there is
no $J_{\theta}$-holomorphic disks for some $\theta$.

In the hyperk\"{a}hler case we have the following.

\begin{lemma}
If $C$\ is a Lagrangian submanifold of a compact hyperk\"{a}hler manifold $M$,
then for all $\theta\in\lbrack0,2\pi)$, with at most one exception,\ there is
no $J_{\theta}$-holomorphic disk in $M$\ bounding $C$.
\end{lemma}

Proof: Suppose $D$\ is a $J_{\theta}$-holomorphic disk in $M$ with $\partial
D\subset C$, say $\theta=0$. For the complex structure $J_{0}=I$, $\omega
_{K}+i\omega_{J}$ is a holomorphic two form on $M$ and therefore restricts to
zero of on any $I$-holomorphic disk. On the other hand $\omega_{I}$ is a
K\"{a}hler form and hence it is positive on $D$. Therefore we have
\[
\omega_{J}=\omega_{K}=0,\omega_{I}>0\text{,}%
\]
on $D$. We consider the integration of $\omega_{I}$ and $\omega_{K}$ on $D$,
since these two forms restrict to zero on $C$ we have well-defined
homomorphisms
\begin{align*}
\int\omega_{I}  &  :H_{2}\left(  M,C;\mathbb{Z}\right)  \rightarrow
\mathbb{R}\text{,}\\
\int\omega_{K}  &  :H_{2}\left(  M,C;\mathbb{Z}\right)  \rightarrow
\mathbb{R}\text{.}%
\end{align*}

From earlier discussions $\left[  D,\partial D\right]  $\ represents a class
in $H_{2}\left(  M,C;\mathbb{Z}\right)  $ which\ must lie in the kernel of
$\int\omega_{K}$ and not in the kernel of $\int\omega_{I}$, in fact $\int
_{D}\omega_{I}>0$. Therefore, with a fix class in $H_{2}\left(  M,C;\mathbb{Z}%
\right)  $, there is at most one $\theta$\ which can support $J_{\theta}%
$-holomorphic disks representing the given class. $\blacksquare$

\bigskip

Remark on the vanishing of the Maslov index: Recall in the real symplectic
geometry, the Maslov index plays a very important role, for instance in
determining the dimension of the space of holomorphic disks. The origin of the
Maslov index is the isomorphism $\pi_{1}\left(  U\left(  n\right)  /O\left(
n\right)  \right)  \cong\mathbb{Z}$, where the space $U\left(  n\right)
/O\left(  n\right)  $ parametrizes linear Lagrangian subspaces in
$\mathbb{R}^{2n}$. For complex linear Lagrangian subspaces in $\mathbb{C}^{2n}
$, their parameter space equals $Sp\left(  n\right)  /U\left(  n\right)
\subset U\left(  2n\right)  /O\left(  2n\right)  $, which has trivial
fundamental group
\[
\pi_{1}\left(  \frac{Sp\left(  n\right)  }{U\left(  n\right)  }\right)
=0\text{.}%
\]
Therefore the Maslov index is always zero in our situation.

\bigskip

\textbf{Definition of Lagrangian category}

\begin{definition}
Given any hyperk\"{a}hler manifold $M,$ we define a category $\mathcal{C}_{M}
$ or simply $\mathcal{C}$, called the Lagrangian category of $M$ as follows:
An object in $\mathcal{C}$ is a Lagrangian subvariety of $M$. Given two
objects $C_{1},C_{2}\in\mathcal{C}^{obj}$ we define the space of morphisms to
be the $\mathbb{Z}$-graded Abelian group
\[
Hom_{\mathcal{C}}\left(  C_{1},C_{2}\right)  ^{\left[  k\right]  }=Ext_{O_{M}%
}^{k}\left(  O_{C_{1}},O_{C_{2}}\right)  \text{.}%
\]
\end{definition}

\bigskip

The composition of morphisms is given by the natural product structure on the
$Ext$'s groups.
\[
Ext_{O_{M}}^{k}\left(  O_{C_{1}},O_{C_{2}}\right)  \otimes Ext_{O_{M}}%
^{l}\left(  O_{C_{2}},O_{C_{3}}\right)  \rightarrow Ext_{O_{M}}^{k+l}\left(
O_{C_{1}},O_{C_{3}}\right)  \text{.}%
\]
Because of the Serre duality and $K_{M}=O_{M}$, the above category carries a
natural \textit{duality} property,
\[
Ext_{O_{M}}^{k}\left(  O_{C_{1}},O_{C_{2}}\right)  \cong Ext_{O_{M}}%
^{2n-k}\left(  O_{C_{2}},O_{C_{1}}\right)  ^{\ast}\text{,}%
\]
provided that $M$ is compact.

\bigskip

\textbf{A symplectic 2-category}

In our situation an object is a Lagrangian $C$ in a fix hyperk\"{a}hler
manifold $M$. In \textit{real }symplectic geometry Weinstein (\cite{BW} and
\cite{We}) defines a symplectic category whose objects are symplectic
manifolds and morphisms are immersed Lagrangian submanifolds inside
$M_{2}\times\bar{M}_{1}$. Here $\bar{M}$ denote the symplectic manifold $M$
with the symplectic form $-\omega$. When $M$ is a hyperk\"{a}hler manifold,
the complex structure of $\bar{M}$ becomes $-J$. In fact we can combine the
two approaches together and define a \textit{symplectic 2-category}: The
objects are hyperk\"{a}hler manifolds, 1-morphisms from $M_{1}$ to $M_{2}$ are
Lagrangian subvarieties in $M_{2}\times\bar{M}_{1}$, 2-morphisms between two
Lagrangian subvarieties $C_{1},C_{2}\subset M_{2}\times\bar{M}_{1}$ are given
by $Ext_{O_{M_{2}\times\bar{M}_{1}}}^{\ast}\left(  O_{C_{1}},O_{C_{2}}\right)
$.

\bigskip

\textbf{The category of Lagrangian coherent sheaves}

The Lagrangian category of $M$ is geometric in nature but it does not have
very good functorial properties. Therefore we also need another category.

\begin{definition}
Let $M$ be a projective hyperk\"{a}hler manifold. Let $D^{b}\left(  M\right)
$ be the derived category of coherent sheaves on $M$. We define the category
of Lagrangian coherent sheaves $D_{Lag}^{b}\left(  M\right)  $ to be the
subcategory of $D^{b}\left(  M\right)  $ generated by those coherent sheaves
$\mathcal{S}$ satisfying $ch\left(  \mathcal{S}\right)  \cup\Omega\in
\oplus_{k>2n+2}H^{k}\left(  M,\mathbb{C}\right)  $
\end{definition}

By the Hard Lefschetz $sl_{2}$-action using $L_{\Omega}$ and $\Lambda_{\Omega
}$, the condition $ch\left(  \mathcal{S}\right)  \cup\Omega\in\oplus
_{k>2n+2}H^{k}\left(  M,\mathbb{C}\right)  $ implies that $ch_{k}\left(
\mathcal{S}\right)  =0$ for $k<n$ and $ch_{n}\left(  \mathcal{S}\right)  $ is
a $\Omega$-primitive cohomology class. In particular the $n$ dimensional
support of $\mathcal{S}$ is a Lagrangian in $M$. For example if $C$ is a
Lagrangian subvariety in $M$ then $O_{C}$ is an object in $D_{Lag}^{b}\left(
M\right)  $.

\subsubsection{Lagrangians intersection}

The Lagrangian category $\mathcal{C}_{M}$ (and similar for $D_{Lag}^{b}\left(
M\right)  $) is closely related to the intersection theory for Lagrangian
subvarieties in $M$.

\begin{theorem}
If $C_{1}$ and $C_{2}$ are two Lagrangian subvarieties of a compact
hyperk\"{a}hler manifold $M$ then
\[%
{\textstyle\sum_{k}}
\dim\left(  -1\right)  ^{k}Ext_{O_{M}}^{k}\left(  O_{C_{1}},O_{C_{2}}\right)
=\left(  -1\right)  ^{n}C_{1}\cdot C_{2}.
\]
\end{theorem}

Proof: We recall the Riemann-Roch formula for the global Ext groups: For any
coherent sheaves $S_{1}$ and $S_{2}$ on $M$ we have,
\[
\dim\left(  -1\right)  ^{k}Ext_{O_{M}}^{k}\left(  S_{1},S_{2}\right)
=\int_{M}\overline{ch}\left(  S_{1}\right)  ch\left(  S_{2}\right)  Td_{M}%
\]
where $\overline{ch}\left(  S_{1}\right)  =\Sigma\left(  -1\right)  ^{k}%
ch_{k}\left(  S_{1}\right)  $. For $S_{i}=O_{C_{i}}$ the structure sheaf of a
subvariety $C_{i}$ of dimension $n$, we have
\begin{align*}
ch_{k}\left(  O_{C_{i}}\right)   &  =0\text{ for }k<n,\\
ch_{n}\left(  O_{C_{i}}\right)   &  =\left[  C_{i}\right]  \text{.}%
\end{align*}
Therefore
\begin{align*}
&  \dim\left(  -1\right)  ^{k}Ext_{O_{M}}^{k}\left(  O_{C_{1}},O_{C_{2}%
}\right) \\
&  =\int_{M}\overline{ch}\left(  O_{C_{1}}\right)  ch\left(  O_{C_{2}}\right)
Td_{M}\\
&  =\int_{M}\left(  \left(  -1\right)  ^{n}\left[  C_{1}\right]
+h.o.t.\right)  \left(  \left[  C_{2}\right]  +h.o.t.\right)  \left(
1+h.o.t.\right) \\
&  =\left(  -1\right)  ^{n}\int_{M}\left[  C_{1}\right]  \cup\left[
C_{2}\right]  =\left(  -1\right)  ^{n}C_{1}\cdot C_{2}\text{.}%
\end{align*}
Here $h.o.t.$ refers to \textit{higher order terms} which do not contribute to
the outcome of the integral. Hence the result. $\blacksquare$

\bigskip

If $C_{1}$ and $C_{2}$ intersect cleanly along $D$ then $C_{1}\cdot C_{2}%
\ $equals the Euler characteristic of $D$ up to sign. To prove this we use the
following useful lemma whose proof is standard.

\begin{lemma}
If $C_{1}$ and $C_{2}$ are two Lagrangian submanifolds of a hyperk\"{a}hler
manifold $M$ which intersect cleanly along $D=C_{1}\cap C_{2}$, then the
symplectic form on $M$ induces a non-degenerate pairing
\[
N_{D/C_{1}}\otimes N_{D/C_{2}}\rightarrow O_{D}\text{.}%
\]
\end{lemma}

Proof: For any $v\in T_{C_{1}}$\ and $w\in T_{D}\subset T_{C_{2}}$ we have
$\Omega\left(  v,w\right)  =0$\ because $D\subset C_{1}$. This implies that
the above pairing is well-defined. Suppose that $v\in T_{C_{1}}$\ is such that
$\Omega\left(  v,u\right)  =0$\ for every $u\in T_{C_{2}}$. Then $v\in\left(
T_{C_{2}}\right)  ^{\bot}=T_{C_{2}}$\ because $C_{2}$\ is Lagrangian. That is
$v\in T_{C_{1}}\cap T_{C_{2}}=T_{D}$. This shows the nondegeneracy of the
pairing. Hence the lemma. $\blacksquare$

\bigskip

In particular we have a natural isomorphism
\[
N_{D/C_{1}}\cong N_{D/C_{2}}^{\ast}\text{.}%
\]
We also have natural exact sequences,
\begin{align*}
0  &  \rightarrow N_{D/C_{2}}\rightarrow\left(  N_{C_{1}/M}\right)
|_{D}\rightarrow T_{D}^{\ast}\rightarrow0,\\
0  &  \rightarrow N_{D/C_{1}}\rightarrow\left(  N_{C_{2}/M}\right)
|_{D}\rightarrow T_{D}^{\ast}\rightarrow0,\\
0  &  \rightarrow N_{D/M}\rightarrow\left(  N_{C_{1}/M}\oplus N_{C_{2}%
/M}\right)  |_{D}\rightarrow T_{D}^{\ast}\rightarrow0,
\end{align*}
comparing relative normal bundles. One should recall that $N_{C_{i}/M}\cong
T_{C_{i}}^{\ast}$ for $i=1,2.$

On the other hand the standard intersection theory gives us that
\begin{align*}
C_{1}\cdot C_{2} &  =e\left(  \left(  N_{C_{1}/M}\oplus N_{C_{2}/M}\right)
|_{D}-N_{D/M}\right) \\
&  =e\left(  T_{D}^{\ast}\right) \\
&  =\left(  -1\right)  ^{\dim D}e\left(  C_{1}\cap C_{2}\right)  \text{.}%
\end{align*}

\begin{corollary}
If $C_{1}$ and $C_{2}$ are two Lagrangian submanifolds in a compact
hyperk\"{a}hler manifold $M$ which intersect cleanly, then
\[
C_{1}\cdot C_{2}=\left(  -1\right)  ^{\dim C_{1}\cap C_{2}}e\left(  C_{1}\cap
C_{2}\right)  \text{.}%
\]
\end{corollary}

Remark: When $M=T^{\ast}X$\ and $C_{i}$ is the conormal bundle $N_{S_{i}%
/X}^{\ast}$ of a submanifold $S_{i}\subset X$. The intersection $C_{1}\cap
C_{2}$\ is compact if and only if it lies inside the zero section $X\subset
M$. This happens precisely when $T_{x}S_{1}+T_{x}S_{2}=T_{x}X$\ for all $x\in
S_{1}\cap S_{2}$, i.e. $S_{1}$\ intersects $S_{2}$\ transversely. In this
situation $C_{1}\cdot C_{2}$\ makes sense and equals to $\left(  -1\right)
^{\dim S_{1}\cap S_{2}}e\left(  S_{1}\cap S_{2}\right)  $.

\bigskip

When $C_{1}$ and $C_{2}$ intersect cleanly we expect that individual
$Ext_{O_{M}}^{k}\left(  O_{C_{1}},O_{C_{2}}\right)  $ can be computed, not
just their Euler characteristics which is given by $\left(  -1\right)
^{n}C_{1}\cdot C_{2}$. For example we have following results.

\begin{theorem}
If $C$ is a Lagrangian submanifold of a compact hyperk\"{a}hler manifold $M$,
then there is an isomorphism of vector spaces,
\[
Ext_{O_{M}}^{k}\left(  O_{C},O_{C}\right)  \cong H^{k}\left(  C,\mathbb{C}%
\right)  ,
\]
for all $k$ provided that the normal bundle of $C$ can be extended to the
whole $M$. For instance it holds true when $C$ is a complete intersection in
$M$.
\end{theorem}

Proof: We consider a Koszul resolution of $O_{C}$ in $M$,
\[
0\rightarrow\Lambda^{n}E\rightarrow\Lambda^{n-1}E\rightarrow\cdots\rightarrow
E\rightarrow O_{M}\rightarrow O_{C}\rightarrow0\text{,}%
\]
where $E$ is a vector bundle on $M$ whose restriction to $C$ is the conormal
bundle $N_{C/M}^{\ast}$ which is isomorphic to the tangent bundle of $C$ by
the Lagrangian condition. The groups $Ext_{O_{M}}^{\ast}\left(  O_{C}%
,O_{C}\right)  $ can be computed as the hypercohomology of the complex of
sheaves \underline{$\mathit{Hom}$}$_{O_{M}}\left(  \Lambda^{\ast}%
E,O_{C}\right)  $ ( \cite{GH}). Note that
\begin{align*}
\underline{\mathit{Hom}}_{O_{M}}\left(  \Lambda^{q}E,O_{C}\right)   &
\cong\underline{\mathit{Hom}}_{O_{M}}\left(  O_{M},\Lambda^{q}E^{\ast}\otimes
O_{C}\right) \\
&  \cong\underline{\mathit{Hom}}_{O_{C}}\left(  O_{C},\Lambda^{q}T_{C}^{\ast
}\right) \\
&  \cong\Omega^{q}\left(  C\right)  \text{.}%
\end{align*}
From the definition of the Koszul complex, the restriction of its dual complex
$\Lambda^{\ast}E^{\ast}$ to $C$ has trivial differentials. Therefore
\begin{align*}
&  Ext_{O_{M}}^{k}\left(  O_{C},O_{C}\right) \\
&  \cong\mathbb{H}\left(  0\rightarrow O_{C}\overset{0}{\rightarrow}\Omega
^{1}\left(  C\right)  \overset{0}{\rightarrow}\cdots\overset{0}{\rightarrow
}\Omega^{n}\left(  C\right)  \rightarrow0\right) \\
&  \cong%
{\textstyle\bigoplus_{p+q=k}}
H^{p,q}\left(  C\right) \\
&  \cong H^{k}\left(  C,\mathbb{C}\right)  \text{.}%
\end{align*}
The last equality uses the K\"{a}hlerian property of $C$. Hence the result.
$\blacksquare$

\bigskip

\begin{theorem}
If $C_{1}$ and $C_{2}$ are two Lagrangian submanifolds of a compact
hyperk\"{a}hler manifold $M$ which intersect transversely along $C_{1}\cap
C_{2}=\left\{  p_{1},...,p_{s}\right\}  $, then there is an isomorphism of
vector spaces,
\[
Ext_{O_{M}}^{k}\left(  O_{C_{1}},O_{C_{2}}\right)  \cong\left\{
\begin{array}
[c]{cccc}%
0 &  & \text{if }k\neq n & \\
&  &  & \\%
{\textstyle\bigoplus_{p_{i}}}
\mathbb{C}\cong\mathbb{C}^{s} &  & \text{if }k=n\text{.} &
\end{array}
\right.
\]
\end{theorem}

Proof: Without loss of generality we can assume that $C_{1}\cap C_{2}=\left\{
p\right\}  .$ We take a Koszul resolution of $O_{C_{1}}$ as before,
\[
0\rightarrow\Lambda^{n}E\rightarrow\Lambda^{n-1}E\rightarrow\cdots\rightarrow
E\rightarrow O_{M}\rightarrow O_{C_{1}}\rightarrow0\text{.}%
\]
Because $C_{1}$ and $C_{2}$ intersect transversely, the restriction of this
resolution to $C_{2}$ gives a Koszul resolution of $O_{p}$ in $C_{2},$%
\[
0\rightarrow\Lambda^{n}F\rightarrow\Lambda^{n-1}F\rightarrow\cdots\rightarrow
F\rightarrow O_{C_{2}}\rightarrow O_{p}\rightarrow0\text{,}%
\]
where $F$ is the restriction of $E$ to $C_{2}$. Now%

\begin{align*}
\mathit{Hom}_{O_{M}}\left(  \Lambda^{q}E,O_{C_{2}}\right)   &  \cong
\mathit{Hom}_{O_{M}}\left(  O_{M},\Lambda^{q}E^{\ast}\otimes O_{C_{2}}\right)
\\
&  \cong\mathit{Hom}_{O_{C_{2}}}\left(  O_{C_{2}},\Lambda^{q}F^{\ast}\right)
\\
&  \cong\mathit{Hom}_{O_{C_{2}}}\left(  \Lambda^{q}F,O_{C_{2}}\right)  .
\end{align*}
Therefore
\begin{align*}
Ext_{O_{M}}^{k}\left(  O_{C_{1}},O_{C_{2}}\right)   &  \cong\mathbb{H}%
^{k}\left(  \mathit{Hom}_{O_{M}}\left(  \Lambda^{\ast}E,O_{C_{2}}\right)
\right) \\
&  \cong\mathbb{H}^{k}\left(  \mathit{Hom}_{O_{C_{2}}}\left(  \Lambda
^{q}F,O_{C_{2}}\right)  \right) \\
&  \cong Ext_{O_{C_{2}}}^{k}\left(  O_{p},O_{C_{2}}\right)  \text{.}%
\end{align*}

The theorem follows from the above equation and basic properties of the
extension groups \cite{GH}. $\blacksquare$

\bigskip

\textbf{Intersection Euler characteristics}

From previous discussions the Euler characteristic of any smooth subvariety
$S$ in $\mathbb{P}^{n}$ equals
\[
\chi\left(  S\right)  =\left(  -1\right)  ^{\dim S}N_{S/\mathbb{P}^{n}}^{\ast
}\cdot\mathbb{P}^{n}\text{.}%
\]
Even though the intersection of the two Lagrangians $N_{S/\mathbb{P}^{n}%
}^{\ast}$ and $\mathbb{P}^{n}$ is taken place in a non-compact manifold,
namely $T^{\ast}\mathbb{P}^{n}$, their intersection number is well-defined in
this case because their intersection occurs inside a compact region and
moreover they do not intersect even \textit{at infinity}. For singular variety
in $\mathbb{P}^{n}$ we use this as a definition.

\begin{definition}
For any subvariety $S$ in $\mathbb{P}^{n}$ we define its intersection Euler
characteristic to be the following intersection number inside $T^{\ast
}\mathbb{P}^{n}$,
\[
\bar{\chi}\left(  S\right)  =\left(  -1\right)  ^{\dim S}N_{S/\mathbb{P}^{n}%
}^{\ast}\cdot\mathbb{P}^{n}\text{.}%
\]
\end{definition}

It is not difficult to see that if $S$ is a smooth simple normal crossing
subvariety of $\mathbb{P}^{n}$ then its intersection Euler characteristic
$\bar{\chi}\left(  S\right)  $ equals to the usual Euler characteristic of its
normalization. For plane curves, such $S$ would have only double point
singularities. By local computations, we can prove the following formula for
intersection Euler characteristics for plane curves which might even have cusp
singularities: For any plane curve $S$ of degree $d$ with $\delta$ doubles
points, $\kappa$ cusps and no other singularities, we have
\[
\bar{\chi}\left(  S\right)  =d^{2}-3d+2\delta+3\kappa\text{.}%
\]

\bigskip

\subsubsection{\label{Sec Red Proj functor}Reduction functor and projection functor}

Reduction of a symplectic manifold $M$ is induced from a coisotropic
submanifold $D$ in $M$. Given any Lagrangian subvariety in $M$ we can
construct another one which lives inside $D$, called the \textit{projection}
and also a Lagrangian subvariety in the reduced symplectic space, called the
\textit{reduction}. We will review these constructions and we will used them
later in the normalized Legendre transformation.

\bigskip

Given any coisotropic submanifold $D$ in $M$ it has an integrable distribution
$\left(  T_{D}\right)  ^{\bot}\subset T_{D}$. We denote the natural projection
to the leave space $B$ as
\[
\pi_{D}:D\rightarrow B.
\]

If we assume $B$ is smooth then it has a natural holomorphic symplectic
structure and $B$ is called the \textit{reduction }of $D$, for simplicity we
assume $B$ is also hyperk\"{a}hler. The existence of such symplectic structure
on $B$ and also our later construction of two transformations it induced are
based on the following linear algebra lemma,

\begin{lemma}
Suppose $M$ is a symplectic vector space with symplectic form $\Omega$, $D$ is
a coisotropic subspace and $C$ is a Lagrangian subspace. Then

(1) $\Omega|_{D}$ induced a symplectic structure on $D/D^{\bot}$;

(2) $C\cap D+D^{\bot}\subset D$ is a Lagrangian subspace in $M$;

(3) $\left(  C\cap D\right)  /\left(  C\cap D^{\bot}\right)  $ is a Lagrangian
subspace in $D/D^{\bot}$.
\end{lemma}

The proof of this lemma is standard and readers can find it in chapter 5 of
\cite{BW} for instance.

\bigskip

\textbf{Reduction and projection of a Lagrangian}

Suppose that $C$ is a Lagrangian subvariety of $M$ and we denote its smooth
locus as $C^{sm}$. We construct a subvariety $C^{red}$ in $B$ (resp.
$C^{proj}$ in $M$) called the \textit{reduction }of $C$ (resp.
\textit{projection }of $C$) as follow,
\[
C^{red}=\overline{\pi_{D}\left(  C^{sm}\cap D\right)  }\subset B.
\]
and
\[%
\begin{array}
[c]{ccccc}%
C^{proj} & \subset &  D & \subset &  M\\
\downarrow & \square & \downarrow &  & \\
C^{red} & \subset &  B. &  &
\end{array}
\]
Because of the above linear algebra lemma, both $C^{red}$ in $B$ and
$C^{proj}$ in $M$ are Lagrangian subvarieties. This has not yet define
functors on Lagrangian categories because it is not so easy to see how to
construct the functor on the morphism level. For this purpose the derived
category of Lagrangian coherent sheaves $D_{Lag}^{b}\left(  \bullet\right)  $
serves a better role.

We define a functor between derived categories called the \textit{reduction
functor} as follow: We consider the subvariety $B\times_{B}D\subset B\times M
$ and denote the projection morphism from $B\times M$ to its first and second
factor as $\pi_{B}$ and $\pi_{M}$ respectively then we define
\begin{align*}
\mathbf{R}_{D}  &  :D_{Lag}^{b}\left(  M\right)  \rightarrow D_{Lag}%
^{b}\left(  B\right) \\
\mathbf{R}_{D}\left(  \bullet\right)   &  =R\pi_{B\ast}\left(  O_{B\times
_{B}D}\overset{L}{\otimes}\pi_{M}^{\ast}\left(  \bullet\right)  \right)
\text{.}%
\end{align*}
\qquad

We can also define a \textit{projection functor }as follow: We consider a
subvariety $D\times_{B}D\subset M\times M$ and denote the projection morphisms
from $M\times M$ to its first and second factors as $\pi_{1}$ and $\pi_{2}$
then we define
\begin{align*}
\mathbf{P}_{D}  &  :D_{Lag}^{b}\left(  M\right)  \rightarrow D_{Lag}%
^{b}\left(  M\right) \\
\mathbf{P}_{D}\left(  \bullet\right)   &  =R\pi_{1\ast}\left(  O_{D\times
_{B}D}\overset{L}{\otimes}\pi_{2}^{\ast}\left(  \bullet\right)  \right)
\text{.}%
\end{align*}
It can be checked that the image of any Lagrangian coherent sheaf under
$\mathbf{P}_{D}$ or $\mathbf{R}_{D}$ is indeed a Lagrangian coherent
sheaf.\footnote{To be precise we should talk about complex of sheaves in the
dervied category.}

Remark: We expect that the functors $\mathbf{R}_{D}$ and $\mathbf{P}_{D}$
induced from a coisotropic submanifold $D$ to enjoy many good properties. For
example $\mathbf{R}_{D}\circ\mathbf{R}_{D}^{\ast}=1$, $\mathbf{P}_{D}%
\circ\mathbf{P}_{D}=\mathbf{P}_{D}$ and $\mathbf{P}_{D}=\mathbf{P}_{D}^{\ast}$
for some suitably defined adjoint functors $\mathbf{R}_{D}^{\ast}$ and
$\mathbf{P}_{D}^{\ast}$. We also expect that $\mathbf{R}_{D} $ is an injective
functor and $\mathbf{P}_{D}$ is a surjective functor. Moreover one should be
able to relax the smoothness assumption on $D$ and $B$.

\section{\label{Sec Legendre}Legendre transformation}

In this section we will study a Legendre transformation of Lagrangian
subvarieties along a coisotropic exceptional subvariety in $M$ and a Plucker
type formula. We will start by recalling the Legendre transformation in the
classical mechanics. It can be reinterpreted as a transformation of conormal
bundles. Then we descend this transformation to the hyperk\"{a}hler quotient
$T^{\ast}\mathbb{P}^{n}$. It has a natural generalization to any Mukai
elementary modification.

We establish a \textit{Plucker type formula} which relates intersection
numbers of Lagrangian subvarieties under the Legendre transformation. Then we
define a \textit{normalized Legendre transformation }which enjoys much better
functorial properties. The definition of this normalized transformation uses
the reduction and the projection of the Lagrangian subvariety with respect to
the coisotropic exceptional locus.

\subsection{Classical Legendre transformation}

The origin of the Legendre transformation is from the classical mechanics. It
is a transformation from the Lagrangian mechanics on $TX$ to the Hamiltonian
mechanics on $T^{\ast}X$. Given a function $L=L\left(  q,v\right)
:TX\rightarrow\mathbb{R}$ usually called a \textit{Lagrangian,} we can form a
closed two form on $TX$ as follows,
\[
\omega_{L}=\sum\frac{\partial^{2}L}{\partial q^{i}\partial v_{j}}dq^{i}\wedge
dq^{j}+\frac{\partial^{2}L}{\partial v_{i}\partial v_{j}}dv_{i}\wedge
dq^{j}\text{.}%
\]
Suppose $L$ is \textit{non-degenerate} in the sense that $\det\left(
\frac{\partial^{2}L}{\partial v_{i}\partial v_{j}}\right)  \neq0$ at every
point of $TX$, the \textit{Legendre transformation} is defined as the
following map,
\begin{align*}
\mathcal{L} &  :TX\rightarrow T^{\ast}X\\
\left(  q,v\right)   &  \rightarrow\left(  q,\frac{\partial L}{\partial v_{i}%
}dq^{i}\right)  .
\end{align*}
It pulls back the canonical symplectic form $\omega$ on $T^{\ast}X$ to the
above $\omega_{L}$ on $TX$, i.e. $\mathcal{L}^{\ast}\omega=\omega_{L}$. We
also define a \textit{Hamiltonian,}
\begin{align*}
H &  :T^{\ast}X\rightarrow\mathbb{R}\\
H &  =%
{\textstyle\sum}
v^{i}\frac{\partial L}{\partial v^{i}}-L\text{,}%
\end{align*}
and it is called the \textit{Legendre transformation} of the function $L$.

\bigskip

\textbf{Legendre transform on vector spaces}

From above discussions, we saw that the Legendre transformation is really a
fiberwise transformation from $T_{q}X$ to $T_{q}^{\ast}X$ using a function
$L\left(  q,\cdot\right)  $ whose Hessian is non-degenerate at every point. We
can reformulate it as follow: Suppose $V$ is a finite dimensional vector space
and $f:V\rightarrow\mathbb{C}$ is a function such that its Hessian is
non-degenerate at every point,
\[
\det\left(  \frac{\partial^{2}f}{\partial x^{i}\partial x^{j}}\right)
\neq0\text{.}%
\]
We look at the graph of $df$ in $T^{\ast}V=V\times V^{\ast}$ then the
\textit{Legendre transformation }induced by $f$ is the map,
\begin{align*}
\mathcal{L}_{f}  &  :V\rightarrow V^{\ast},\\
\mathcal{L}_{f}  &  =\pi_{V^{\ast}}\circ df,
\end{align*}
where $\pi_{V^{\ast}}$ is the projection from $V\times V^{\ast}$ to its second
factor $V^{\ast}$. That is $\mathcal{L}_{f}\left(  x\right)  =\xi$ if and only
if $\xi_{i}=\frac{\partial f}{\partial x^{i}}$ for all $i$ in local coordinates.

We also define the Legendre transformation of the function $f$ on $V$ as a
function
\begin{align*}
f^{\vee}  &  :V^{\ast}\rightarrow\mathbb{C},\\
f^{\vee}\left(  \xi\right)   &  =\Sigma x^{i}\xi_{i}-f\left(  x\right)
\text{,}%
\end{align*}
with $\xi=\mathcal{L}_{f}\left(  x\right)  $.

\bigskip

Geometrically this transformation arises from the natural isomorphism,
\[
T^{\ast}V\cong T^{\ast}\left(  V^{\ast}\right)  ,
\]
for any finite dimensional vector space $V$. The reason is $\left(  V^{\ast
}\right)  ^{\ast}=V$ and therefore both sides are isomorphic to $V\times
V^{\ast}$.

Now the graph $C$ of $df$ in $T^{\ast}V$ is a Lagrangian submanifold of
$T^{\ast}V$ with its canonical symplectic form. Using the above isomorphism we
can treat
\[
C\subset T^{\ast}\left(  V^{\ast}\right)  \text{,}%
\]
as another Lagrangian submanifold.\footnote{Under the natural isomorphism
$T^{\ast}V\cong T^{\ast}\left(  V^{\ast}\right)  $, their canonical symplectic
forms are identified up to a minus sign. In particular they have identical
Lagrangian submanifolds.} Under the non-degeneracy assumption $C$ is also a
graph of a function on $V^{\ast}$, this function is precisely the above
$f^{\vee}$.

It is obvious that the Legendre transformation is involutive, i.e.
\begin{align*}
\mathcal{L}_{f}\circ\mathcal{L}_{f^{\vee}}  &  =1_{V}\\
\mathcal{L}_{f^{\vee}}\circ\mathcal{L}_{f}  &  =1_{V^{\ast}}\\
\left(  f^{\vee}\right)  ^{\vee}  &  =f\text{.}%
\end{align*}

\bigskip

Remark: There is another point of view for the Legendre transformation,
explained in Guillemin's book \cite{Gu} as follow: Suppose, on a real vector
space, $f:V\rightarrow\mathbb{R}$ is a strictly convex function with a
critical point, which is necessarily unique and a global minimum, such a $f$
is called \textit{stable.} Then $\mathcal{L}_{f}\left(  x\right)  =\xi$ if and
only if $x$ is the unique critical point for the function $f\left(  x\right)
-\xi\left(  x\right)  $. Moreover $f^{\vee}\left(  \xi\right)  =-\min_{x\in
V}f_{\xi}\left(  x\right)  $. We can also identify the image of $\mathcal{L}%
_{f}$ as $\left\{  \xi\in V^{\ast}:f\left(  x\right)  -\xi\left(  x\right)
\text{ is stable}\right\}  $.

\bigskip

\textbf{A different perspective: conormal bundles}

Recall that we associate to any function $f$ on $V$ a Lagrangian submanifold
in $T^{\ast}V$ given by the graph of $df$, then we use the natural isomorphism
$T^{\ast}V\cong T^{\ast}\left(  V^{\ast}\right)  $ to define the Legendre
transformation $f^{\vee}$. Besides the graph of $df$, there is another natural
Lagrangian submanifold in $T^{\ast}V$ associated to $f$, namely the conormal
bundle of the zero set of $f$,
\[
N_{S/V}^{\ast}\subset T^{\ast}V,
\]
where $S=\left\{  x\in V:f\left(  x\right)  =0\right\}  $. When $f$ is
\textit{homogeneous} this Lagrangian submanifold, considered inside $T^{\ast
}\left(  V^{\ast}\right)  $, turns out to be the conormal bundle of the zero
set of $f^{\vee}$. This allows us to define the Legendre transformation in
much greater generality, at least in the homogeneous case which is just the
right setting for the projective geometry.

\begin{theorem}
\label{Thm Transf conormal bdl}Suppose $f:V\rightarrow\mathbb{C}$ is a
homogenous polynomial such that its zero set $S=\left\{  x\in V:f\left(
x\right)  =0\right\}  $ is smooth. Then its conormal bundle $N_{S/V}^{\ast
}\subset T^{\ast}V$, when viewed as a submanifold of $T^{\ast}\left(  V^{\ast
}\right)  $, equals $N_{S^{\vee}/V^{\ast}}^{\ast}$ the conormal variety of
$S^{\vee}=\left\{  \xi\in V^{\ast}:f^{\vee}\left(  \xi\right)  =0\right\}
\subset V^{\ast}$.
\end{theorem}

Proof: When $f$ is a homogeneous polynomial of degree $p$, we have
\[
\Sigma x^{i}\frac{\partial f}{\partial x^{i}}=pf\left(  x\right)  \text{.}%
\]
Under the Legendre transformation $\xi_{i}=\frac{\partial f}{\partial x^{i}%
}\left(  x\right)  $ this gives,
\begin{align*}
f^{\vee}\left(  \xi\right)   &  =\Sigma x^{i}\xi_{i}-f\left(  x\right) \\
&  =pf\left(  x\right)  -f\left(  x\right) \\
&  =\left(  p-1\right)  f\left(  x\right)  \text{.}%
\end{align*}

The conormal bundle of $S$ is given by,%

\[
N_{S/V}^{\ast}=\left\{  \left(  x,\eta\right)  \in V\times V^{\ast}:f\left(
x\right)  =0\text{ and }\eta_{i}=c\frac{\partial f}{\partial x^{i}}\left(
x\right)  \text{ for all }i\text{, for some }c\right\}  .
\]
So we need to verify that the same set can be described as
\[
\left\{  \left(  x,\eta\right)  \in V\times V^{\ast}:f^{\vee}\left(
\eta\right)  =0\text{ and }x^{i}=b\frac{\partial f^{\vee}}{\partial\xi_{i}%
}\left(  \eta\right)  \text{ for all }i\text{, for some }b\right\}  .
\]

Now we suppose $\left(  x,\eta\right)  \in N_{S/V}^{\ast}$. Since $f$ is
homogenous of degree $p,$ we have
\[
\frac{\partial f}{\partial x^{i}}\left(  ex\right)  =e^{p-1}\frac{\partial
f}{\partial x^{i}}\left(  x\right)  ,
\]
for any number $e$. Therefore $\eta_{i}=c\frac{\partial f}{\partial x^{i}%
}\left(  x\right)  =\frac{\partial f}{\partial x^{i}}\left(  c^{\prime
}x\right)  $ for some constant $c^{\prime}$. That is
\[
\eta=L_{f}\left(  c^{\prime}x\right)  \text{.}%
\]
Hence
\[
f^{\vee}\left(  \eta\right)  =f^{\vee}\left(  L_{f}\left(  c^{\prime}x\right)
\right)  =\left(  p-1\right)  f\left(  c^{\prime}x\right)  =\left(
p-1\right)  \left(  c^{\prime}\right)  ^{p}f\left(  x\right)  =0\text{.}%
\]
Similarly using the inverse Legendre transformation, we have
\[
\frac{\partial f^{\vee}}{\partial\xi_{i}}\left(  \eta\right)  =c^{\prime}%
x^{i}\text{.}%
\]
Hence the result. $\blacksquare$

\bigskip

One can also find an indirect proof of this in \cite{GKZ}. Because of this
result, we can now define the Legendre transformation for any finite set of
homogenous polynomials via the conormal bundle of their common zero set. This
approach works particularly well for projective spaces and it is closely
related to the dual variety construction.

\bigskip

\subsection{Legendre transform in hyperk\"{a}hler manifolds}

\subsubsection{Legendre transform in $T^{\ast}\mathbb{P}^{n}$ and dual varieties}

Recall that the hyperk\"{a}hler structure on $T^{\ast}\mathbb{P}^{n}$ can be
constructed as the hyperk\"{a}hler quotient of $T^{\ast}V$ with $V=\mathbb{C}%
^{n+1}$ by the natural $S^{1}$-action,
\[
T^{\ast}\mathbb{P}^{n}=\left\{  \left(  x,\xi\right)  \in V\times V^{\ast}%
:\xi\left(  x\right)  =1,\left|  x\right|  ^{2}-\left|  \xi\right|
^{2}=1\right\}  /S^{1}\text{.}%
\]
The Legendre transformation on a linear symplectic space $T^{\ast}V$ comes
from the natural isomorphism,
\[
T^{\ast}V\cong T^{\ast}\left(  V^{\ast}\right)  \text{.}%
\]
We are going to descend this transformation to their hyperk\"{a}hler quotients
$T^{\ast}\mathbb{P}^{n}$ and $T^{\ast}\mathbb{P}^{n\ast}$. In particular we
only need to look at those functions on $V$ which are homogenous.

The first issue is the natural isomorphism between $T^{\ast}V$ and $T^{\ast
}\left(  V^{\ast}\right)  $ does not descend to their hyperk\"{a}hler
quotients.\ Instead we have a natural \textit{birational }map between
$T^{\ast}\mathbb{P}^{n}$ and $T^{\ast}\mathbb{P}^{n\ast}$ which preserves
their symplectic structures. To see this we recall that
\[
T^{\ast}\mathbb{P}^{n\ast}=\left\{  \left(  \xi,x\right)  \in V^{\ast}\times
V:\xi\left(  x\right)  =1,\left|  \xi\right|  ^{2}-\left|  x\right|
^{2}=1\right\}  /S^{1},
\]
and the birational map is given by
\begin{align*}
\Phi &  :T^{\ast}\mathbb{P}^{n}\dashrightarrow T^{\ast}\mathbb{P}^{n\ast}\\
\Phi\left(  x,\xi\right)   &  =\left(  \left|  \frac{x}{\xi}\right|
\xi,\left|  \frac{\xi}{x}\right|  x\right)  .
\end{align*}
It can be verified directly $\Phi$ is an isomorphism outside their zero
sections, which are given by $\xi=0$ in $T^{\ast}\mathbb{P}^{n}$ and $x=0$ in
$T^{\ast}\mathbb{P}^{n\ast}$. Moreover $\Phi$ pullbacks the canonical
symplectic structure on $T^{\ast}\mathbb{P}^{n\ast}$ to the one on $T^{\ast
}\mathbb{P}^{n}$ outside their zero sections.

\bigskip

Because of the theorem \ref{Thm Transf conormal bdl} we define the
\textit{Legendre transformation} in the projective setting as follow: For any
homogenous function $f:V\rightarrow\mathbb{C}$ which defines a smooth
hypersurface $S=\left\{  f=0\right\}  \subset\mathbb{P}^{n}$, the function
itself can be recovered from its conormal bundle $N_{S/\mathbb{P}^{n}}^{\ast}$
inside $T^{\ast}\mathbb{P}^{n}$. The Legendre transform $f^{\vee}:V^{\ast
}\rightarrow\mathbb{C}$ defines the dual hypersurface $S^{\vee}\subset
\mathbb{P}^{n\ast}$ under the non-degenerate assumption and we have
\[
N_{S^{\vee}/\mathbb{P}^{n\ast}}=\overline{\Phi\left(  N_{S/\mathbb{P}^{n}%
}\backslash\mathbb{P}^{n}\right)  }\text{,}%
\]
where $\mathbb{P}^{n}$ denote the zero section in $T^{\ast}\mathbb{P}^{n}$.

An arbitrary Lagrangian subvariety $C$ in $T^{\ast}\mathbb{P}^{n}$ can be
regarded as a generalized homogenous function on $V$ unless $C$ is the zero
section $\mathbb{P}^{n}$. Motivated from above discussions, we define the
\textit{Legendre transformation} $C^{\vee}$ of $C$ as follow,
\[
C^{\vee}=\overline{\Phi\left(  C\backslash\mathbb{P}^{n}\right)  \text{.}}%
\]

It has the following immediate properties: (i) $C^{\vee}$ is a Lagrangian
subvariety of $T^{\ast}\mathbb{P}^{n\ast}$; (ii) $C\subset T^{\ast}%
\mathbb{P}^{n}$ and $C^{\vee}\subset T^{\ast}\mathbb{P}^{n\ast}$ are
isomorphic outside the zero sections; (iii) the inversion property $\left(
C^{\vee}\right)  ^{\vee}=C$.

\bigskip

Remark: Recall that $T^{\ast}\mathbb{P}^{n}$ is the hyperk\"{a}hler quotient
of $T^{\ast}V$ by the natural $S^{1}$ action. In fact the Legendre
transformation on $T^{\ast}\mathbb{P}^{n}$ can be regarded as a $S^{1}%
$-invariant Legendre transformation on $T^{\ast}V$ using the symplectic
quotient by $S^{1}$: If $C$ is a Lagrangian in $T^{\ast}\mathbb{P}^{n}$ then
there is a unique $\mathbb{C}^{\times}$-invariant Lagrangian subvariety of
$T^{\ast}V,$ denote $D$ such that $C$ is the symplectic quotient of $D$ by
$S^{1}$. The Legendre transformation of $D$ will be a Lagrangian subvariety
$D^{\vee}$ in $T^{\ast}\left(  V^{\ast}\right)  $ which is again
$\mathbb{C}^{\times}$-invariant. The symplectic quotient of $D^{\vee}$ by
$S^{1}$ would be our transformation $C^{\vee}$ in $T^{\ast}\mathbb{P}^{n\ast}$.

\bigskip

\textbf{Dual varieties in }$\mathbb{P}^{n}$ \textbf{and a Plucker formula}

For any subvariety $S$ in $\mathbb{P}^{n}$ we can associated a dual variety
$S^{\vee}$ in the dual projective space $\mathbb{P}^{n\ast}$. The dual variety
$S^{\vee}$ is the closure of all hyperplanes in $\mathbb{P}^{n}$ which are
tangent to some smooth point in $S$. From our previous discussions, the
conormal variety of $S^{\vee}$ is simply the Legendre transformation of the
conormal variety of $S$. In this sense our Legendre transformation is a
generalization of the dual variety construction. The relationship between dual
varieties and the Legendre transformation has been briefly addressed before in
various places.

As an example the dual variety of the Fermat hypersurface
\[
S=\left\{  x_{0}^{p}=x_{1}^{p}+\cdots+x_{n}^{p}\right\}  \subset\mathbb{P}^{n}%
\]
is
\[
S^{\vee}=\left\{  \xi_{0}^{q}=\xi_{1}^{q}+\cdots+\xi_{n}^{q}\right\}
\subset\mathbb{P}^{n\ast}%
\]
where $\frac{1}{p}+\frac{1}{q}=1$. Notice that $q$ is no longer an integer if
$p>2$, in fact $S^{\vee}$ is a hypersurface of degree $p\left(  p-1\right)
^{n-1}$.

The study of dual varieties is a very rich classical subject in algebraic
geometry (see \cite{GKZ} for instance). We list a few known facts (1) the
inversion property $\left(  S^{\vee}\right)  ^{\vee}=S$, (2) $S$ is
irreducible if and only if $S^{\vee}$ is irreducible, (3) If $x\in S$ and
$\xi\in S^{\vee}$ are smooth points, then $\xi$ is tangent to $S$ at $x$ if
and only if $x$ is tangent to $S^{\vee}$ at $\xi$.

For plane curves, there are Pl\"{u}cker formulae which related various
geometric quantities between $S$ and $S^{\vee}$: Suppose $S\subset
\mathbb{P}^{2}$ is a plane curve of degree $d$ with $\delta$ double points,
$\kappa$ cusps and no other singularities, we denote the corresponding
quantities for $S^{\vee}$ as $d^{\vee},\delta^{\vee},\kappa^{\vee}$. Then
Pl\"{u}cker formulae say,
\begin{align*}
d^{\vee} &  =d\left(  d-1\right)  -2\delta-3\kappa\\
\kappa^{\vee} &  =3d^{2}-6d-6\delta-8\kappa\text{.}%
\end{align*}
A similar formula in the higher dimensional setting is obtained by Kleiman in
\cite{Kl}.

In the next section we will discuss a similar formula for the Legendre
transformation in any hyperk\"{a}hler manifold. In the case of conormal
varieties inside $T^{\ast}\mathbb{P}^{n}$ the formula says: For any
subvarieties $S_{i}\subset\mathbb{P}^{n}$ of dimension $s_{i}$ we denote their
conormal varieties as $C_{i}\subset T^{\ast}\mathbb{P}^{n}$ then we have
\[
C_{1}\cdot C_{2}+\frac{\left(  C_{1}\cdot\mathbb{P}^{n}\right)  \left(
C_{2}\cdot\mathbb{P}^{n}\right)  }{\left(  -1\right)  ^{n+1}\left(
n+1\right)  }=C_{1}^{\vee}\cdot C_{2}^{\vee}+\frac{\left(  C_{1}^{\vee}%
\cdot\mathbb{P}^{n\ast}\right)  \left(  C_{2}^{\vee}\cdot\mathbb{P}^{n\ast
}\right)  }{\left(  -1\right)  ^{n+1}\left(  n+1\right)  },
\]
or in terms of the intersection Euler characteristics $\bar{\chi}$ we have
\[
C_{1}\cdot C_{2}\pm\frac{1}{n+1}\bar{\chi}\left(  S_{1}\right)  \bar{\chi
}\left(  S_{2}\right)  =C_{1}^{\vee}\cdot C_{2}^{\vee}\pm\frac{1}{n+1}%
\bar{\chi}\left(  S_{1}^{\vee}\right)  \bar{\chi}\left(  S_{2}^{\vee}\right)
.
\]
provided that $S_{1}$ and $S_{2}$ intersect transversely and the same for
their duals.

In the special case of plane curves, the intersection number of their conormal
varieties $C_{1}\cdot C_{2}$ is simply the product of the degree of the curves
$d_{1}d_{2}$. In fact the proof of the formula can be reduced to the case when
$S_{2}$ is a point and it reads as follow: For any plane curve $S$ we have%

\[
3d^{\vee}=-\bar{\chi}\left(  S\right)  -2\bar{\chi}\left(  S^{\vee}\right)  .
\]
When $S$ has only double point and cusp singularities, the above formula can
be proven by the Pl\"{u}cker formulae. Conversely the Pl\"{u}cker formula
$d^{\vee}=d\left(  d-1\right)  -2\delta-3\kappa$ also follows from it and our
earlier formula for $\bar{\chi}\left(  S\right)  $.

\subsubsection{Legendre transform along $\mathbb{P}^{n}$ and a Pl\"{u}cker
type formula}

Now we study the Legendre transformation on a general hyperk\"{a}hler manifold
$M$ of dimension $2n\geq4$.

\bigskip

\textbf{Flop along }$\mathbb{P}^{n}$

Recall that every embedded $\mathbb{P}^{n}$ in $M$ is a Lagrangian submanifold
and it has a neighborhood $U$ which is symplectomorphic to a neighborhood of
the zero section in $T^{\ast}\mathbb{P}^{n}$ and we continue to denote it by
$U$. Therefore we can \textit{flop} such a $\mathbb{P}^{n}$ in $M$ to obtain
another holomorphic symplectic manifold $M^{\prime} $. To see how this surgery
work, we look at the birational transformation $\Phi:T^{\ast}\mathbb{P}%
^{n}\dashrightarrow T^{\ast}\mathbb{P}^{n\ast}$ and let $U^{\prime}$ be the
image of $U$, i.e. $U^{\prime}=\overline{\Phi\left(  U\right)  }$. Then
$M^{\prime}=\left(  M\backslash U\right)  \cup U^{\prime}$.\footnote{To be
precise with the holomorphic structure on $M^{\prime}$ we should write
$M^{\prime}=\left(  M\backslash\overline{U_{0}}\right)  \cup U^{\prime}$ for
some open set $U_{0}\supset\mathbb{P}^{n}$ satisfying $\overline{U_{0}}\subset
U$.} Since $M$ and $M^{\prime}$ are isomorphic outside a codimension $n$
subspace, $M^{\prime}$ inherits a holomorphic two form $\Omega^{\prime}$ from
$M$ by the Hartog's theorem. Moreover being a section of the canonical line
bundle and non-vanishing outside a codimension $n$ subset, $\left(
\Omega^{\prime}\right)  ^{n}$ must be non-vanishing everywhere. That is
$M^{\prime}$ is a holomorphic symplectic manifold with the symplectic form
$\Omega^{\prime}$. In most cases $M^{\prime}$ is actually a hyperk\"{a}hler
manifold, nevertheless Namikawa had a simple example \cite{Na} showing that
this cannot always be the case.

We will denote the natural birational map between $M$ and $M^{\prime}$ as
\[
\Phi_{M}:M\dashrightarrow M^{\prime}.
\]

\bigskip

\textbf{Examples of flop}

(1) (Mukai) Suppose $X$ is a degree two K3 surface, i.e. $X$ is a double cover
of $\mathbb{P}^{2}$ branched along a sextic curve, $\pi:X\rightarrow
\mathbb{P}^{2}$. We denote $J_{0}$ the degree zero compactified Picard scheme
for degree 2 curves in $X$. There is a birational map
\[
\Phi:S^{\left[  2\right]  }X\dashrightarrow J_{0}%
\]
which associates to any 2 generic points $\left(  p_{1},p_{2}\right)  $ in $X$
to the unique degree two curve $C$ passing through them together with the line
bundle $\omega_{C}\otimes O\left(  -p_{1}-p_{2}\right)  $.

Note that $\mathbb{P}^{2}$ embeds inside $S^{\left[  2\right]  }X$ via the
preimage of $\pi$. In fact $\Phi$ is the flop of $S^{\left[  2\right]  }X$
along this $\mathbb{P}^{2}$.

(2) (Beauville): Let $X\subset\mathbb{P}^{3}$ be a smooth quartic surface. Any
2 points on $X$ defines a line in $\mathbb{P}^{3}$ which intersects $X$ at two
other points. This defines a birational map
\[
\Phi:S^{\left[  2\right]  }X\dashrightarrow S^{\left[  2\right]  }X.
\]
If $X$ does not contain any line, this is an isomorphism. If $X$ contains $k$
lines $L_{1},...,L_{k}$, then $\Phi$ is the flop of $S^{\left[  2\right]  }X$
along $S^{\left[  2\right]  }L_{j}\cong\mathbb{P}^{2}$'s.

\bigskip

\textbf{Legendre transformation and a Pl\"{u}cker type formula}

Suppose $M$ and $M^{\prime}$ are hyperk\"{a}hler manifolds which are related
by a flop $\Phi_{M}:M\dashrightarrow M^{\prime}$ along a $\mathbb{P}%
^{n}\subset M$. For any Lagrangian subvariety $C$ in $M$ which does not
contain $\mathbb{P}^{n}$ we define its Legendre transform to be the Lagrangian
subvariety $C^{\vee}$ in $M^{\prime}$ defined as,
\[
C^{\vee}=\overline{\Phi_{M}\left(  C\backslash\mathbb{P}^{n}\right)  }.
\]
Since $\left(  C_{1}\cup C_{2}\right)  ^{\vee}=C_{1}^{\vee}\cup C_{2}^{\vee}
$, we can extend the definition of the Legendre transformation to the free
Abelian group generated by all Lagrangian subvarieties of $M$ except
$\mathbb{P}^{n}$.

Clearly the Legendre transformation has the inversion property, namely
$\left(  C^{\vee}\right)  ^{\vee}=C$. However the Legendre transformation does
not preserve the intersection numbers, i.e. $C_{1}\cdot C_{2}\neq C_{1}^{\vee
}\cdot C_{2}^{\vee}$. Instead they satisfy the following Pl\"{u}cker type
formula \cite{Le}.

\begin{theorem}
Suppose $\Phi_{M}:M\dashrightarrow M^{\prime}$ is a flop along $\mathbb{P}%
^{n}$ between projective hyperk\"{a}hler manifolds then
\[
C_{1}\cdot C_{2}+\frac{\left(  C_{1}\cdot\mathbb{P}^{n}\right)  \left(
C_{2}\cdot\mathbb{P}^{n}\right)  }{\left(  -1\right)  ^{n+1}\left(
n+1\right)  }=C_{1}^{\vee}\cdot C_{2}^{\vee}+\frac{\left(  C_{1}^{\vee}%
\cdot\mathbb{P}^{n\ast}\right)  \left(  C_{2}^{\vee}\cdot\mathbb{P}^{n\ast
}\right)  }{\left(  -1\right)  ^{n+1}\left(  n+1\right)  }\text{,}%
\]
for any Lagrangian subvarieties $C_{1}$ and $C_{2}$ not containing
$\mathbb{P}^{n}$.
\end{theorem}

In the above Pl\"{u}cker type formula, it is more natural to interpret the LHS
(and similar for the RHS) as follow,
\begin{align*}
&  C_{1}\cdot C_{2}+\frac{\left(  C_{1}\cdot\mathbb{P}^{n}\right)  \left(
C_{2}\cdot\mathbb{P}^{n}\right)  }{\left(  -1\right)  ^{n+1}\left(
n+1\right)  }\\
&  =\left(  C_{1}-\frac{\left(  C_{1}\cdot\mathbb{P}^{n}\right)  }{\left(
\mathbb{P}^{n}\cdot\mathbb{P}^{n}\right)  }\mathbb{P}^{n}\right)  \cdot\left(
C_{2}-\frac{\left(  C_{2}\cdot\mathbb{P}^{n}\right)  }{\left(  \mathbb{P}%
^{n}\cdot\mathbb{P}^{n}\right)  }\mathbb{P}^{n}\right)  \text{.}%
\end{align*}
Note that $\left(  C-\frac{\left(  C\cdot\mathbb{P}^{n}\right)  }{\left(
\mathbb{P}^{n}\cdot\mathbb{P}^{n}\right)  }\mathbb{P}^{n}\right)
\cdot\mathbb{P}^{n}=0$ for any Lagrangian subvariety $C$ in $M$, including
$C=\mathbb{P}^{n}$. Roughly speaking the LHS (resp. RHS) of the Pl\"{u}cker
type formula is the intersection number of two Lagrangian subvarieties which
do not intersect $\mathbb{P}^{n}$, and therefore the Legendre transformation
along such a $\mathbb{P}^{n}$ should have no effect to their intersection
numbers, thus giving the Pl\"{u}cker type formula in a heuristic way.

\bigskip

\textbf{Normalized Legendre transformation}

It is suggested from above discussions that we should modify the
transformation so that
\[
C-\frac{\left(  C\cdot\mathbb{P}^{n}\right)  }{\left(  \mathbb{P}^{n}%
\cdot\mathbb{P}^{n}\right)  }\mathbb{P}^{n}\rightarrow C^{\vee}-\frac{\left(
C^{\vee}\cdot\mathbb{P}^{n\ast}\right)  }{\left(  \mathbb{P}^{n\ast}%
\cdot\mathbb{P}^{n\ast}\right)  }\mathbb{P}^{n\ast}.
\]
We also want transform $\mathbb{P}^{n}$, the center of the flop. Thus we
arrive to the following definition of a \textit{normalized Legendre
transformation }$\mathcal{L}$:
\begin{align*}
\mathcal{L}\left(  C\right)   &  =C^{\vee}+\frac{\left(  C\cdot\mathbb{P}%
^{n}\right)  +\left(  -1\right)  ^{n+1}\left(  C^{\vee}\cdot\mathbb{P}^{n\ast
}\right)  }{n+1}\mathbb{P}^{n\ast}\text{ if }C\neq\mathbb{P}^{n}\\
\mathcal{L}\left(  \mathbb{P}^{n}\right)   &  =\left(  -1\right)
^{n}\mathbb{P}^{n\ast}\text{.}%
\end{align*}
Now our Pl\"{u}cker type formula can be rephrased as the following simple
identity,
\[
C_{1}\cdot C_{2}=\mathcal{L}\left(  C_{1}\right)  \cdot\mathcal{L}\left(
C_{2}\right)  \text{.}%
\]

\bigskip

Unlike our earlier Legendre transformation, this normalized Legendre
transformation is rather non-trivial even when $n=1$. When $M$ is a K3 surface
every embedded $\mathbb{P}^{1}$ is also called an $\left(  -2\right)  $-curve
because $\mathbb{P}^{1}\cdot\mathbb{P}^{1}=-2$. Flopping $\mathbb{P}^{1}$ in
$M$ is trivial, i.e. $M=M^{\prime}$. This is because a point in $\mathbb{P}%
^{1}$ is also a hyperplane. Therefore the Legendre transformation
$C\rightarrow C^{\vee}=C$ is just the identity transformation. However the
normalized Legendre transformation is given by%

\[
\mathcal{L}\left(  C\right)  =C-\left(  C\cdot\mathbb{P}^{1}\right)
\mathbb{P}^{1},
\]
which is well-defined for any cohomology class of the K3 surface $M$ and
induces an automorphism of $H^{2}\left(  M,\mathbb{Z}\right)  $, namely the
reflection with respect to the class $\left[  \mathbb{P}^{1}\right]  $. We can
identify $\mathcal{L}\left(  C\right)  $ as the Dehn twist of $C$ along the
Lagrangian $\mathbb{P}^{1}$, or $S^{2}$, in $M$. If we use all the $\left(
-2\right)  $-curves in the whole twistor family then their corresponding Dehn
twists generate $Aut\left(  H^{2}\left(  M,\mathbb{Z}\right)  \right)  $,
which gives all diffeomorphisms of $M$ up to isotopy by a result of Donaldson.

Remark: It would be interesting to compare our normalized Legendre
transformation with the Dehn twist along Lagrangian $S^{n}$ in the real
symplectic geometry as studied by P. Seidal \cite{Se}. On the level of derived
category, Thomas suggests that $\mathcal{L}$ could be the mirror object to a
Dehn twisted operation defined by him and Seidal \cite{ST}.

\bigskip

Suppose $M$ is a projective hyperk\"{a}hler manifold of dimension $2n$ and
\[
\pi:M\rightarrow Z
\]
is a projective contraction with normal exceptional locus $D$ in $M$. We
assume that $\pi\left(  D\right)  $ is a single point, i.e. $Z$ has isolated
singularity. Then (i) $D\cong\mathbb{P}^{n}$ if $n\geq2$ and (ii) $D$ is an
ADE configuration of $\mathbb{P}^{1}$ if $n=1$. Conversely every such $D$ can
be contracted inside $M$. Moreover the normalized Legendre transformation is
defined in all these cases and very interesting.

\bigskip

\textbf{A categorical transformation}

Now we have a transformation $\mathcal{L}$ which takes Lagrangians in $M$ to
Lagrangians in $M^{\prime}$ and it respects their intersection numbers. It is
then natural to wonder if we have a categorical transformation between the
Lagrangian categories of $M$ and $M^{\prime}$, namely we want to have an
isomorphism on the cohomology groups,
\[
Ext_{O_{M}}^{q}\left(  O_{C_{1}},O_{C_{2}}\right)  \overset{\cong}%
{\rightarrow}Ext_{O_{M^{\prime}}}^{q}\left(  O_{\mathcal{L}\left(
C_{1}\right)  },O_{\mathcal{L}\left(  C_{2}\right)  }\right)
\]
rather than just their Euler characteristics, i.e. intersection numbers.

Suppose $\Phi:M\dashrightarrow M^{\prime}$ is a flop along $P=\mathbb{P}^{n}$
as before. Let $\widetilde{M}$ be the blow up of $M$ along $P$ and we denote
its exceptional divisor as $\widetilde{P}$. $\widetilde{P}$ admits two
$\mathbb{P}^{n-1}$-fibration over $\mathbb{P}^{n}$ provided that $n\geq2.$ A
point in $\widetilde{P}$ is a pair $\left(  p,H\right)  $ with $H\subset
\mathbb{P}^{n}$ a hyperplane and $p\in H$. The two $\mathbb{P}^{n-1}%
$-fibrations correspond to sending the above point to $p$ and $H$
respectively. We can blow down $\widetilde{P}$ along the second fibration to
obtain $M^{\prime}$, we write these morphisms as
\[
M\overset{\pi}{\leftarrow}\widetilde{M}\overset{\pi^{\prime}}{\rightarrow
}M^{\prime}\text{.}%
\]

Now we define a Legendre functor
\begin{align*}
\mathbf{L} &  :D_{Lag}^{b}\left(  M\right)  \rightarrow D_{Lag}^{b}\left(
M^{\prime}\right) \\
\mathbf{L}\left(  \bullet\right)   &  =\pi_{\ast}^{\prime}\pi^{\ast}\left(
\bullet\right)  \text{.}%
\end{align*}
It is not difficult to show that the image of any Lagrangian coherent sheaf on
$M$ under $\mathbf{L}$ is a Lagrangian coherent sheaf on $M^{\prime}$. Using
the machinery developed by Bondal and Orlov in \cite{BO} $\mathbf{L}$ defines
a equivalence of categories \cite{Le}. When $\mathcal{S}=O_{C}$ for a
Lagrangian subvariety $C$ in $M$, then $\mathbf{L}\left(  \mathcal{S}\right)
$ should be closely related to $O_{\mathcal{L}\left(  C\right)  }$. This
transformation will play an important role in the proof of the Pl\"{u}cker
type formula \cite{Le}.

\subsubsection{Legendre transform on Mukai elementary modification}

Now we consider a general projective contraction on $M$, as discussed in
section \ref{Sec Contract Birat}, with smooth discriminant locus $B$,
\[%
\begin{array}
[c]{ccl}%
D & \subset &  M\\
\downarrow &  & \downarrow\pi\\
B & \subset &  Z.
\end{array}
\]
The exceptional locus $D$ is a $\mathbb{P}^{k}$-bundle over $B$, which has a
natural symplectic form. Moreover the normal bundle of $D$ in $M$ is the
relative cotangent bundle of $D\rightarrow B$. The Mukai elementary
modification produces another symplectic manifold $M^{\prime}$ by replacing
the $\mathbb{P}^{k}$-bundle $D\rightarrow B$ with the dual $\mathbb{P}^{k\ast
}$-bundle over $B$.
\[%
\begin{array}
[c]{rcccccl}%
D & \subset &  M & \overset{\Phi_{M}}{\dashrightarrow} & M^{\prime} & \supset
&  D^{\prime}\\
\mathbb{P}^{k}\downarrow &  & \,\downarrow\pi &  & \pi^{\prime}\downarrow\, &
& \downarrow\mathbb{P}^{k\ast}\\
B & \subset &  Z & = & Z & \supset &  B
\end{array}
\]

\bigskip

We define the Legendre transformation of a Lagrangian subvariety $C$ in $M$
not lying inside $D$ as follow,
\[
C^{\vee}=\overline{\Phi\left(  C\backslash D\right)  }\text{.}%
\]

In order to write down the Pl\"{u}cker type formula in this general situation,
we need to recall from section \ref{Sec Red Proj functor} that $C$ determines
(i) a Lagrangian subvariety $C^{proj}$ in $M$ lying inside $D$, called the
\textit{projection} and (ii) a Lagrangian subvariety $C^{red}$ in $B$, called
the \textit{reduction. }The restriction of $\pi$ to $C^{proj}$ is a
$\mathbb{P}^{k}$-bundle over $C^{red}\subset B$ and $C^{proj}=\pi^{-1}\left(
C^{red}\right)  \subset M$. It is not difficult to see that we have $\left(
C^{\vee}\right)  ^{proj}=\left(  \pi^{\prime}\right)  ^{-1}\left(
C^{red}\right)  \subset M^{\prime}$ and
\[
C^{proj}\cdot C^{proj}=C^{\vee proj}\cdot C^{\vee proj}\text{.}%
\]

The \textit{Pl\"{u}cker type formula }in this general case reads as follow,%

\begin{align*}
& \left(  C_{1}-\frac{C_{1}\cdot C_{1}^{proj}}{C_{1}^{proj}\cdot C_{1}^{proj}%
}C_{1}^{proj}\right)  \cdot\left(  C_{2}-\frac{C_{2}\cdot C_{2}^{proj}}%
{C_{2}^{proj}\cdot C_{2}^{proj}}C_{2}^{proj}\right) \\
& =\left(  C_{1}^{\vee}-\frac{C_{1}^{\vee}\cdot C_{1}^{\vee proj}}{C_{1}^{\vee
proj}\cdot C_{1}^{\vee proj}}C_{1}^{proj}\right)  \cdot\left(  C_{2}^{\vee
}-\frac{C_{2}^{\vee}\cdot C_{2}^{\vee proj}}{C_{2}^{\vee proj}\cdot
C_{2}^{\vee proj}}C_{2}^{\vee proj}\right)  .
\end{align*}
It can be proven using the same method as in the previous situation.
\textit{The normalized Legendre transformation }$\mathcal{L}$ for this general
case is defined as follow,
\begin{align*}
\mathcal{L}\left(  C\right)   & =C^{\vee}+\frac{\left(  -1\right)  ^{k}C\cdot
C^{proj}-C^{\vee}\cdot C^{\vee proj}}{C^{proj}\cdot C^{proj}}C^{\vee
proj}\text{ when }C\nsubseteq D\text{,}\\
& =\left(  -1\right)  ^{k}\left(  \pi^{\prime}\right)  ^{-1}\left(
C^{red}\right)  \text{ when }C\subset D\text{.}%
\end{align*}
It preserves intersection products of Lagrangian subvarieties in $M$ and
$M^{\prime}$,
\[
C_{1}\cdot C_{2}=\mathcal{L}\left(  C_{1}\right)  \cdot\mathcal{L}\left(
C_{2}\right)
\]

\section{Appendix}

In this appendix we will recall some facts about hyperk\"{a}hler geometry that
we used in this article. We assume $M$ is a compact hyperk\"{a}hler manifold.

\bigskip

\textbf{Hard Lefschetz }$sl_{2}$\textbf{-action using }$\Omega$

We consider the homomorphism
\[%
\begin{array}
[c]{cc}%
L_{\Omega}: & \Omega^{p,q}\left(  M,\mathbb{C}\right)  \rightarrow
\Omega^{p+2,q}\left(  M,\mathbb{C}\right) \\
& L_{\Omega}\left(  \phi\right)  =\phi\cup\Omega\text{,}%
\end{array}
\]
and its adjoint homomorphism
\[
\Lambda_{\Omega}:\Omega^{p+2,q}\left(  M,\mathbb{C}\right)  \rightarrow
\Omega^{p,q}\left(  M,\mathbb{C}\right)  .
\]
As studied by Fujiki in \cite{Fu} they define a $sl_{2}\left(  \mathbb{C}%
\right)  $ action on $\Omega^{\ast,\ast}\left(  M,\mathbb{C}\right)  $.
Moreover this action can be descended to the cohomology groups $H^{\ast,\ast
}\left(  M,\mathbb{C}\right)  $ because $\Omega$ is a parallel form. We call
this the Hard Lefschetz $sl_{2}$-action on $M$ using $\Omega$. A cohomology
class $\phi\in H^{\ast,\ast}\left(  M,\mathbb{C}\right)  $ is called $\Omega
$-primitive if $\Lambda_{\Omega}\phi=0$. As in the standard Hodge theory for
K\"{a}hler manifolds, we have an $\Omega$-primitive decomposition of the
cohomology groups of $M$. In particular $L_{\Omega}$ in injective on
$H^{p,q}\left(  M,\mathbb{C}\right)  $ provided that $p<n$. This $sl_{2}%
$-action forms part of the $so\left(  4,1\right)  $ action defined by Verbitsky.

\bigskip

\textbf{Bogomolov-Beauville quadratic form}

We normalize $\Omega\ $with $\int_{M}\Omega^{n}\bar{\Omega}^{n}=1$. The
Bogomolov-Beauville (unnormalized) quadratic form $q$ on $H^{2}\left(
M,\mathbb{R}\right)  $ is defined as follows (see \cite{Be} for more
details):
\[
q\left(  \phi\right)  =\frac{n}{2}\int\Omega^{n-1}\bar{\Omega}^{n-1}\phi
^{2}+\left(  1-n\right)  \int\Omega^{n-1}\bar{\Omega}^{n}\phi\int\Omega
^{n}\bar{\Omega}^{n-1}\phi.
\]
It is nondegenerate and has signature $\left(  3,b_{2}-3\right)  $. It is not
difficult to check that if $\beta_{1}$ is an ample class and $\beta_{2}$ is an
effective divisor class then $q\left(  \beta_{1},\beta_{2}\right)  >0.$ Using
Torelli theorem for hyperk\"{a}hler manifolds, one can obtain the following
result of Fujiki \cite{Fu}: If $\alpha\in H^{4k}\left(  M\right)  $ is a
polynomial in Chern classes of $M$, then there exists a constant $c_{\alpha}$
such that
\[
\int_{M}\alpha\phi^{2n-2k}=c_{\alpha}q\left(  \phi\right)  ^{n-k}\text{.}%
\]

\bigskip

\textit{Acknowledgments: The author thanks I.M.S. in the Chinese University of
Hong Kong in providing support and an excellent research environment where
much of this work was carried out, special thank to B. Hassett and A. Todorov
for many useful discussions in I.M.S.. The project is also partially supported
by NSF/DMS-0103355. The author also thank D. Abramovich, T. Bridgeland, R.
Thomas and A. Voronov for very helpful discussions.}

\bigskip

{\small Address: School of Mathematics, University of Minnesota, Minneapolis,
MN 55455, U.S.A..}

\bigskip

\begin{thebibliography}{BHL}
\bibitem[BHL]{BHL}D. Burns, Y. Hu, T. Luo, \emph{HyperKahler Manifolds and
Birational Transformations in dimension 4}, math.AG/0004154

\bibitem[BW]{BW}S. Bates, A. Weinstein, \emph{Lectures on the Geometry of
quantization}, Berkeley Math. Lecture Notes, vol. 8., Amer. Math. Soc. 2000.

\bibitem[Be]{Be}A. Beauville, \emph{Varieties Kahleriennes dont la premiere
classe de Chern est nulle}, J. Diff. Geom. 18 (1983) 755-782.

\bibitem[Bo]{Bo}F. Bogomolov, \emph{Hamiltonian K\"{a}hler manifolds,} Soviet
Math. Dokl. 19 (1978) 1462-1465.

\bibitem[BO]{BO}A. Bondal, D. Orlov, \emph{Semiorthogonal decompositions for
algebraic varieties}, alg-geom/9506012.

\bibitem[BL1]{BL1}J. Bryan, N.C. Leung, \emph{The enumerative geometry of K3
surfaces and modular form, }J.A.M.S., 13 (2000), no. 2, 371-410.

\bibitem[BL2]{BL2}J. Bryan and N.C. Leung,\textit{\ }\emph{Generating
functions for the number of curves on abelian surfaces,} Duke Math. J. 99
(1999), no. 2, 311--328.

\bibitem[BL3]{BL3}J. Bryan and N.C. Leung, \emph{Counting curves on rational
surfaces}, Surveys in Differential Geometry, edited by S.T. Yau (1999) 313-340.

\bibitem[Ca]{Ca}E. Calabi, \emph{Metriques Kahleriennes et fibres
holomorphes}, Ann. Ec. Norm. Sup. 12 (1979) 269-294.

\bibitem[Ca2]{Ca2}E. Calabi, \emph{Isometric families of K\"{a}hler
structures, }The Chern Symposium 1979, Springer-Verlag (1980), 23-40.

\bibitem[DM]{DM}R. Donagi, E. Markman, \emph{Spectral covers, algebraically
completely integrable, Hamiltonian systems, and moduli of bundles}, Integrable
systems and quantum groups (Montecatini Terme, 1993), 1--119, Lecture Notes in
Math., 1620, Springer, Berlin, 1996.

\bibitem[DM2]{DM2}R. Donagi, E. Markman, \emph{Cubics, integrable systems, and
Calabi-Yau threefolds,} Proceedings of the Hirzebruch 65 Conference on
Algebraic Geometry (Ramat Gan, 1993), 199--221, Israel Math. Conf. Proc., 9,
Bar-Ilan Univ., Ramat Gan, 1996.

\bibitem[Fr]{Fr}D. Freed, \emph{Special K\"{a}hler manifolds}, C.M.P. 203
(1999) 31-52.

\bibitem[Fu1]{Fu1}A. Fujiki, \emph{On primitive symplectic compact K\"{a}hler
V-manifolds of dimension four,} in 'Classification of algebraic and analytic
manifolds', K. Ueno (ed.), Progress in Math., Birkhauser 39 (1983) 71-125.

\bibitem[Fu]{Fu}A. Fujiki, \emph{On the de Rham cohomology group of a compact
K\"{a}hler symplectic manifold}, Adv. Stud. Pure Math. 10 (1987), Alg. Geom.,
Sendai, 1985, 105-165.

\bibitem[GKZ]{GKZ}I.M. Gelfand, M. Kapranov, A. Zelevinsky,
\emph{Discriminants, resultants and multidimensional determinants},
Birkh\"{a}user (1994).

\bibitem[GH]{GH}P. Griffiths, J. Harris, \emph{Principles of algebraic
geometry}, Wiley-Interscience (1978).

\bibitem[Gu]{Gu}V. Guillemin, \emph{Moment maps and combinatorial invariants
of Hamiltonian }$T^{n}$\emph{-spaces}, Progress in Mathematics, Vol. 122,
Birkh\"{a}user (1994).

\bibitem[Hi]{Hi}N. Hitchin, \emph{Hyperkahler manifolds}, seminaire Bourbaki,
Asterisque 206 (1992) 137-166.

\bibitem[Hi2]{Hi2}N. Hitchin, \emph{The moduli space of complex Lagrangian
submanifolds}, Sir Michael Atiyah: a great mathematician of the twentieth
century. Asian J. Math. 3 (1999), no. 1, 77--91.

\bibitem[HY]{HY}Y. Hu, S.T. Yau, \emph{Hyperkahler manifolds and birational
transformations}, preprint.

\bibitem[Hu]{Hu}D. Huybrechts,\emph{\ Compact hyperk\"{a}hler manifolds: basic
results.} Invent. math. 135 (1999), 63-113.

\bibitem[Kl]{Kl}S. Kleiman \emph{A generalized Teissier-Pl\"{u}cker formula}.
Classification of algebraic varieties (L'Aquila, 1992), 249--260, Contemp.
Math., 162, Amer. Math. Soc., Providence, RI, 1994.

\bibitem[Le]{Le}N.C. Leung, \emph{A general Pl\"{u}cker formula}, in preparation.

\bibitem[Mar]{Mar}E. Markman, \emph{Brill-Neother duality for moduli spaces of
sheaves on K3 surfaces}, J. Algebraic Geom. 10 (2001), no. 4, 623--694.

\bibitem[Mc]{Mc}R.C. McLean, \emph{Deformation of Calibrated Submanifolds},
Commun. Analy. Geom. \textbf{6} (1998) 705-747.

\bibitem[Mu1]{Mu1}S. Mukai, \emph{Symplectic structure of the moduli space of
sheaves on an abelian or K3 surface, }Invent. Math. 77 (1984) 101-116.

\bibitem[Na]{Na}Y. Namikawa, \emph{Deformation theory of singular symplectic
n-folds}, . Math. Ann. 319 (2001), no. 3, 597--623.

\bibitem[Se]{Se}P. Seidal, \emph{Lagrangian two-spheres can be symplectically
knotted}, J. Differential Geom. 52 (1999), no. 1, 145--171.

\bibitem[ST]{ST}P. Seidel, R. Thomas, \emph{Braid group actions on derived
categories of coherent sheaves}, 108 (2001), no. 1, Duke Math. J., 37--108.

\bibitem[Sh]{Sh}N. Shepherd-Barron, \emph{Long extremal rays and symplectic
resolutions}, preprint.

\bibitem[We]{We}A. Weinstein,\textit{\ }\emph{Lectures on symplectic
manifolds}, Regional conference series in mathematics 29, Amer. Math. Soc.,
Providence, 1977.

\bibitem[Wi]{Wi}J. Wierzba, \emph{Contractions of symplectic varieties,} math.AG/9910130.

\bibitem[Ya]{Ya}S.T. Yau, \emph{On the Ricci curvature of a compact K\"{a}hler
manifold and the complex Monge-Ampere equation I}, Comm. Pure and Appl. Math.
31 (1978) 339-411.

\bibitem[YZ]{YZ}S.T. Yau, E. Zaslow, \emph{BPS states, string duality and
nodal curves on K3}, Nucl. Phys. B, 471 (3), (1996) 503-512.
\end{thebibliography}
\end{document}